\documentclass[a4paper,11pt]{article}
\usepackage{amsmath}
\usepackage{amssymb}

\usepackage{texdraw}
\usepackage{mathrsfs}
\usepackage{epsfig}
\usepackage{amsfonts}

\usepackage[dvips]{psfrag}
\usepackage{graphicx}
\usepackage{epstopdf}
\usepackage{subfigure}
\usepackage{color}
\usepackage{booktabs}
\newcommand{\tabincell}[2]{\begin{tabular}{@{}#1@{}}#2\end{tabular}}

\usepackage{geometry}
\usepackage{hyperref}
\hypersetup{colorlinks=true,
            linkcolor=blue,
            anchorcolor=blue,
            citecolor=blue}

\newtheorem{Definition}{Definition}[section]
\newtheorem{Proposition}{Proposition}[section]
\newtheorem{Lemma}{Lemma}[section]
\newtheorem{Theorem}{Theorem}[section]
\newtheorem{Corollary}{Corollary}[section]
\newtheorem{Remark}{Remark}[section]
\newtheorem{Example}{Example}[section]
\newcommand{\bpf}{{\bf Proof:\ \ }}
\newcommand{\epf}{\mbox{}\hfill $\Box$}

\usepackage{caption}
\usepackage {indentfirst}

\begin{document}
\setlength{\baselineskip}{15.2pt} \pagestyle{myheadings}

\title{\bf The number of limit cycles of Josephson equation}

\author{\small Xiangqin Yu$^{a}$, Hebai Chen$^{b}$, Changjian Liu{$^{a,*}$} \\
{\small $^{a}$ School of Mathematics (Zhuhai), Sun Yat-sen University,}\\
{ \small  Zhuhai, 519082,  P.R. China }\\
{\small $^{b}$ School of Mathematics and Statistics, HNP-LAMA, Central South University, }\\
{\small \sl Changsha, 410083, P. R. China}}

\renewcommand{\thefootnote}{}
\footnotetext{$^*$Corresponding author.
\par
E-mail addresses: yuxq25@mail2.sysu.edu.cn(X. Yu), chen\underline{ }hebai@csu.edu.cn(H. Chen),
\par
liuchangj@mail.sysu.edu.cn(C. Liu).}

\date{}
\maketitle
\begin{abstract} In this paper, the existence and number of non-contractible limit cycles of the Josephson equation $\beta \frac{d^{2}\Phi}{dt^{2}}+(1+\gamma \cos \Phi)\frac{d\Phi}{dt}+\sin \Phi=\alpha$ are studied, where $\phi\in \mathbb S^{1}$ and $(\alpha,\beta,\gamma)\in \mathbb R^{3}$. Concretely, by using some appropriate transformations, we prove that such type of limit cycles are changed to limit cycles of some Abel equation. By developing the methods on limit cycles of Abel equation, we prove that there are at most two non-contractible limit cycles, and the upper bound is sharp. At last, combining with the results of the paper (Chen and Tang, J. Differential Equations, 2020), we show that the sum of the number of contractible and non-contractible limit cycles of the Josephson equation is also at most two, and give the possible configurations of limit cycles when two limit cycles appear.
\end{abstract}


\vspace{1mm}
\section{Introduction and main results\label{intro}}
 One of the core problems in the qualitative theory of ordinary differential equations is the study of limit cycles, which is not only closely related to the second half of Hilbert's $16$th problem \cite{cherkas1976number,li2003hilbert}, but also can be used to describe the actual oscillation phenomena of many models \cite{gasull2016number,gasull2020chebyshev,harko2014travelling}.

In this paper, we consider the Josephson equation as follows
\begin{equation}\label{eq1}
 \beta \frac{d^{2}\Phi}{dt^{2}}+(1+\gamma \cos \Phi)\frac{d\Phi}{dt}+\sin \Phi=\alpha,
\end{equation}
which is a physical model related to voltage sources, where $\phi\in \mathbb  S^{1}$ and $(\alpha,\beta,\gamma)\in \mathbb R^{3}$. For more details, we refer the reader to papers \cite{kautz1980proposed,soulen1978impedance}. When $\beta>0$, for the sake of simplicity,  write $(z,\frac{1}{\sigma^{2}})$ instead of $(\frac{d\Phi}{dt}, \beta)$, and let
$z=\sigma y, \alpha=\bar{\alpha}\sigma, t=\frac{\tau}{\sigma},$
 then system \eqref{eq1} is changed to the following Li\'{e}nard system
\begin{equation}\label{eq2}
\begin{cases}
\frac{d\Phi}{d \tau}=y,\\
\frac{dy}{d \tau}=-\sin \Phi+\sigma[\bar{\alpha}-(1+\gamma \cos \Phi)y].\\
\end{cases}
\end{equation}
For sufficiently small $\sigma>0$, Sanders and Cushman \cite{sanders1986limit} obtained the  bifurcation diagram and phase portraits of system \eqref{eq2} by using the averaging method.

Nevertheless, in practical applications $\sigma$ can be not small, therefore, it is worth studying system \eqref{eq2} for generic $\sigma$. By taking the change of variables $\Phi=x, \bar{\alpha}\sigma=a, \sigma=b, \gamma \sigma=c$, system \eqref{eq2} can be written as
\begin{equation}\label{eq3}
\begin{cases}
\frac{dx}{d \tau}=y,\\
\frac{dy}{d \tau}=-(\sin x-a)-(b+c \cos x)y.\\
\end{cases}
\end{equation}
 Chen and Tang \cite{chen2020global} further analyzes the complete dynamics and complex bifurcation phenomena of system \eqref{eq3}. They firstly show that
  one can only consider the case $(a,b,c)\in \Omega:=[0,+\infty)\times(0,+\infty)\times \mathbb R$, and then gave the following  main results.
\begin{Theorem}\label{mainth1}
When $(a,b,c)\in \Omega$, the following statements hold:
\begin{enumerate}
\item[(i)] System \eqref{eq3} includes the following local and global bifurcation surfaces and curves, and cross-sections of the bifurcation diagram are as shown in Fig.~\ref{Fig.1} for an arbitrarily fixed constant $c_{0}$:
\begin{figure}[!htbp]
\centering
\setcounter{subfigure}{0}
\subfigure[when $c=c_{0}<0$]
{\psfrag{1}{{$o$}}
\psfrag{2}{\textcolor{red}{$BT$}}
\psfrag{3}{$a$}
\psfrag{4}{$S_{1}$}
\psfrag{5}{\textcolor{red}{$SC_{2}$}}
\psfrag{6}{$S_{2}$}
\psfrag{7}{\textcolor{red}{$HL$}}
\psfrag{8}{$S_{3}$}
\psfrag{9}{\textcolor{red}{$SC_{1}$}}
\psfrag{0}{\textcolor{red}{$H$}}
\psfrag{a}{$P_{2}$}
\psfrag{b}{$S_{6}$}
\psfrag{c}{\textcolor{red}{$SN$}}
\psfrag{d}{$S_{7}$}
\psfrag{e}{$P_{1}$}
\psfrag{f}{\textcolor{red}{$HE$}}
\psfrag{g}{$S_{4}$}
\psfrag{h}{$S_{5}$}
\psfrag{l}{$b$}
\psfig{file=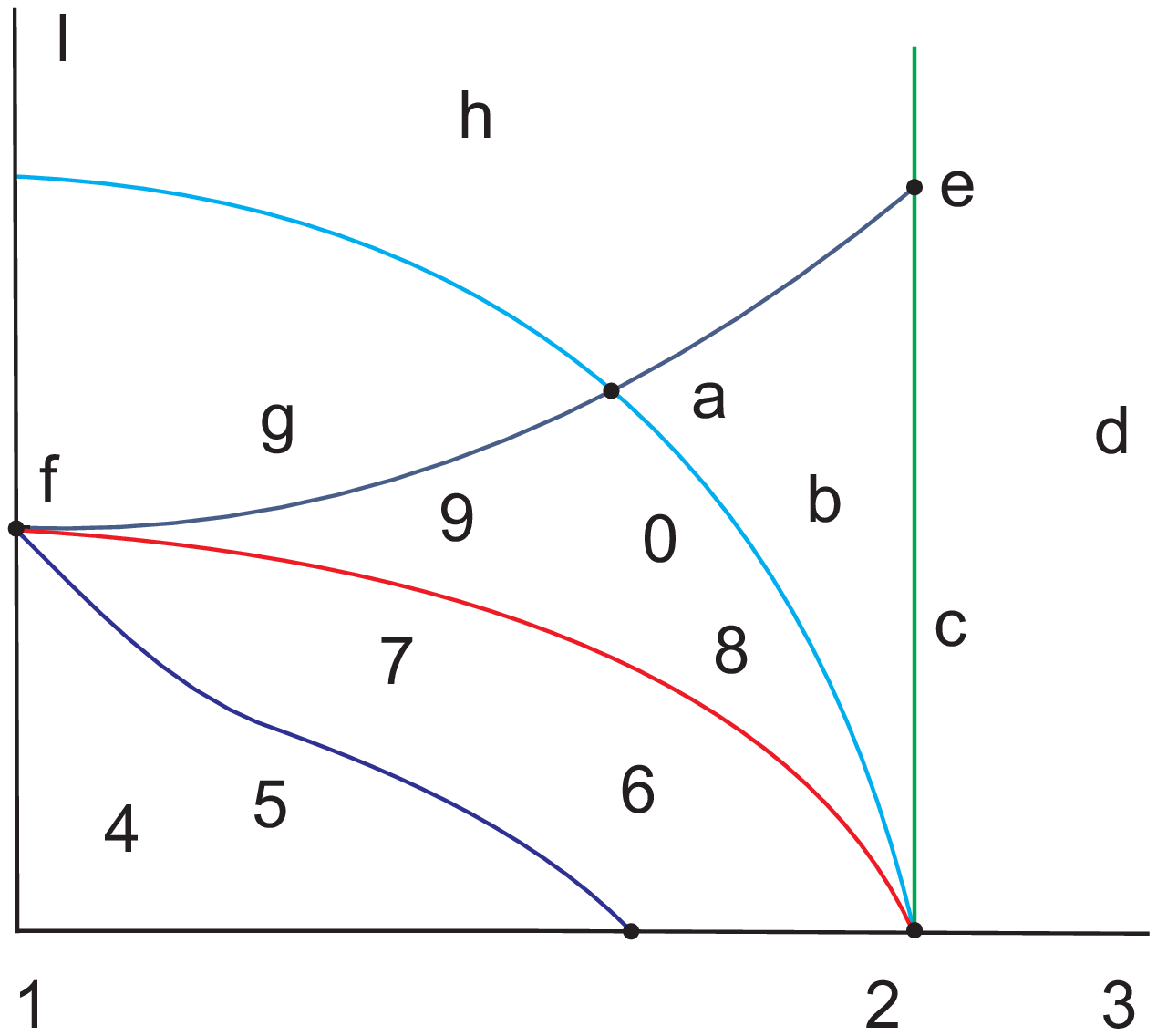,height=1.9in,width=2.4in}}
\hspace{21mm}
\centering
\subfigure[when $c=0$]
{\psfrag{1}{\textcolor{red}{$HEC$}}
\psfrag{7}{\textcolor{red}{$HLC$}}
\psfrag{2}{$1$}
\psfrag{3}{$a$}
\psfrag{4}{$S_{6}$}
\psfrag{5}{\textcolor{red}{$SC_{1}$}}
\psfrag{6}{$S_{5}$}
\psfrag{8}{$S_{7}$}
\psfrag{9}{$P_{1}$}
\psfrag{0}{\textcolor{red}{$SN$}}
\psfrag{a}{$b$}
\psfig{file=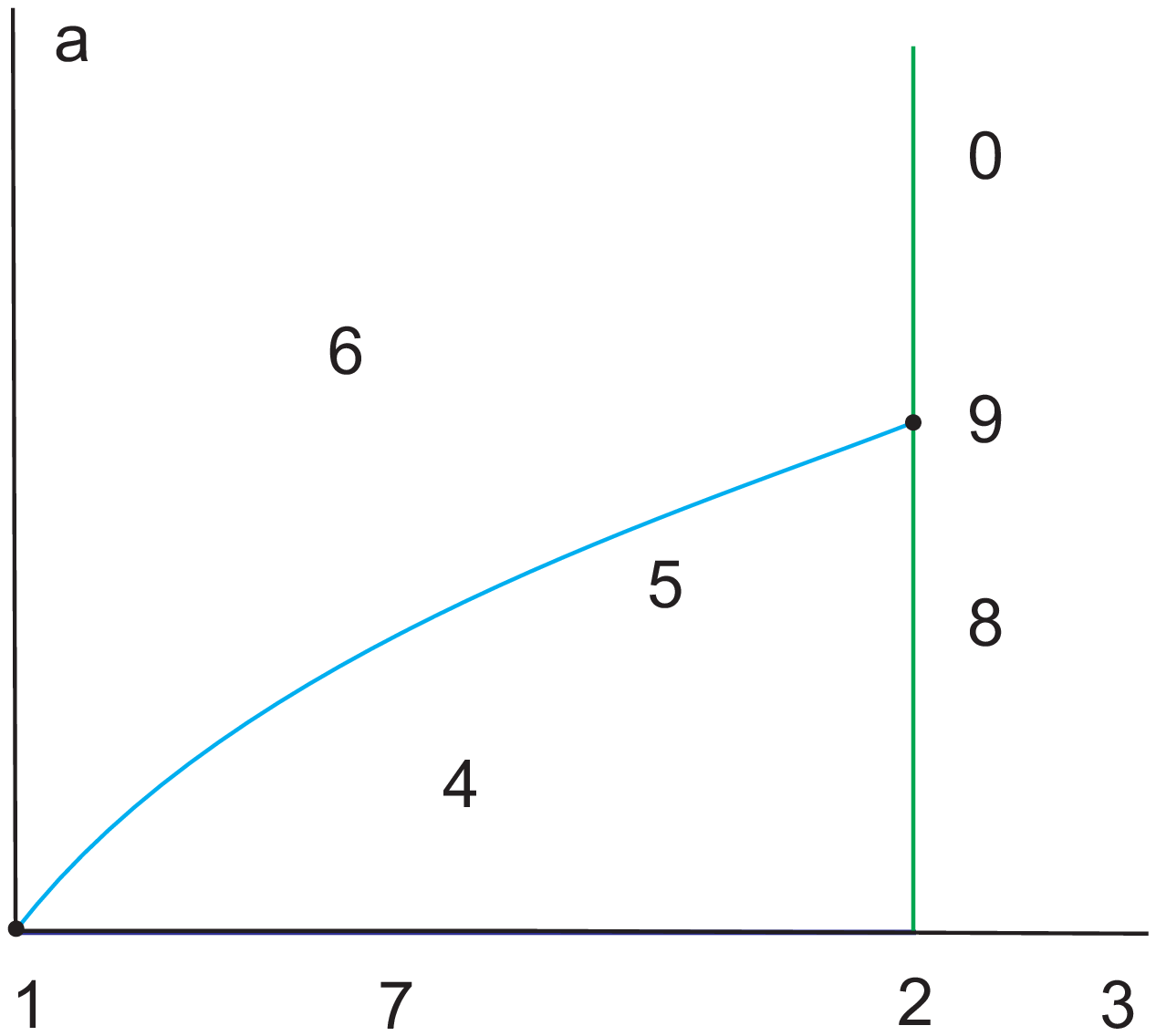,height=1.9in,width=2.4in}}
\\
\centering
\subfigure[when $c=c_{0}>0$]
{\psfrag{1}{\textcolor{red}{$0$}}
\psfrag{2}{$1$}
\psfrag{3}{$a$}
\psfrag{4}{$S_{6}$}
\psfrag{5}{\textcolor{red}{$SC_{1}$}}
\psfrag{6}{$S_{5}$}
\psfrag{7}{$b=1/c$}
\psfrag{8}{$S_{7}$}
\psfrag{9}{$P_{1}$}
\psfrag{0}{\textcolor{red}{$SN$}}
\psfrag{a}{$b$}
\psfig{file=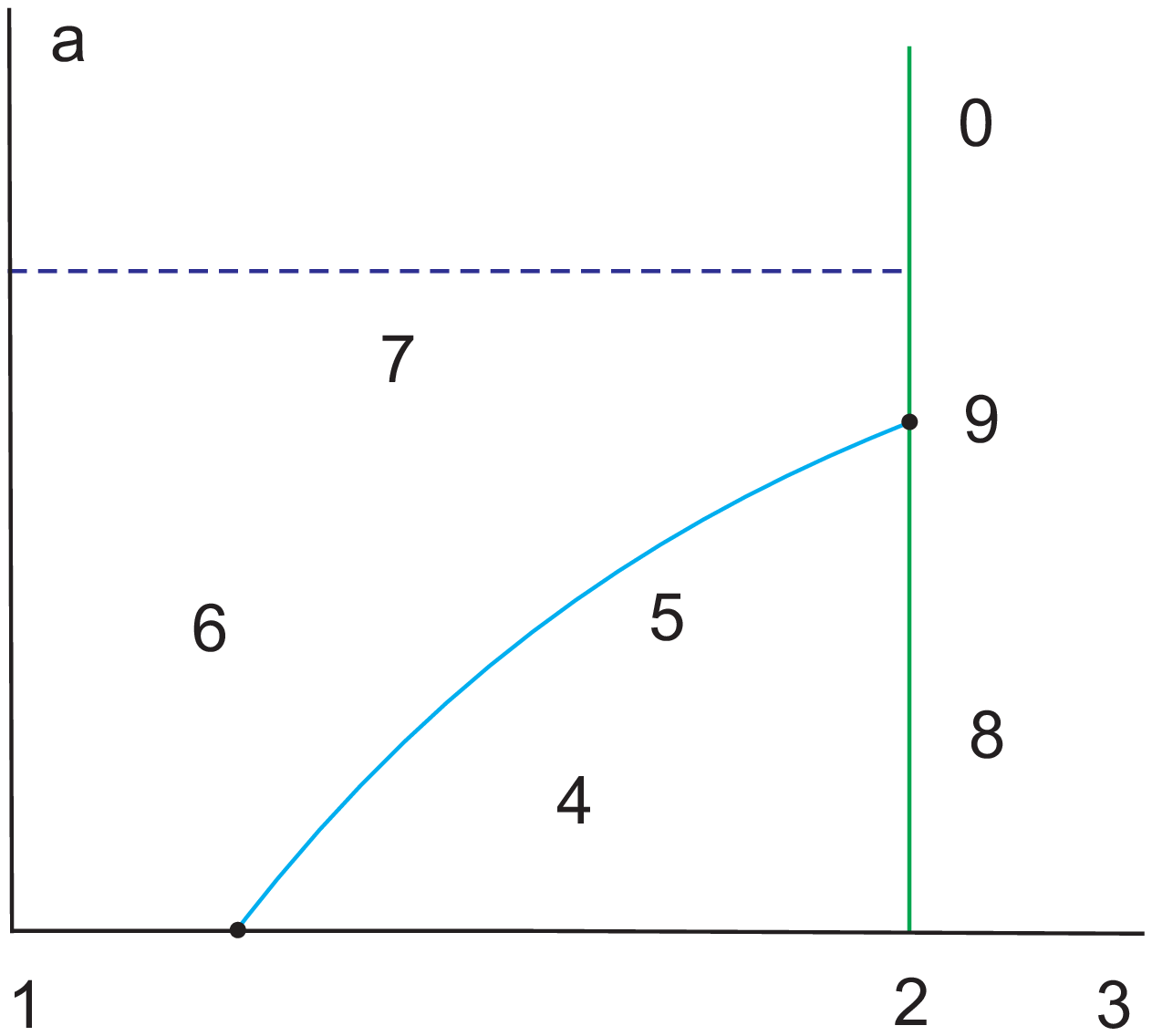,height=1.9in,width=2.4in}}
\begin{center}
\caption{Cross-sections of the bifurcation diagram of system \eqref{eq3}.}
\label{Fig.1}
\end{center}
\end{figure}
\begin{enumerate}
\item[$(a)$] Saddle-node bifurcation surface $$SN:=\{(a,b,c)\in \Omega:a=1\}.$$
\item[$(b)$] Hopf bifurcation surface $$H:=\{(a,b,c)\in \Omega:b=-c\sqrt{1-a^{2}},0\leq a<1\}.$$
\item[$(c)$] Bogdanov-Takens bifurcation curve $$BT:=\{(a,b,c)\in \Omega:a=1,b=0,c<0\}.$$
\item[$(d)$] The bifurcation surface of the homoclinic loop $$HL:=\{(a,b,c)\in \Omega:b=\varphi(a,c),0<a<1,c<0\},$$
where $\varphi\in C^{0},\varphi(1,c)=0,0<\varphi(a,c)<-c\sqrt{1-a^{2}}$ for $0<a<1$ and $c<0$, and $$\varphi(a,c)=\frac{-5c\sqrt{2(1-a)}}{7}+o(-c\sqrt{1-a})$$  when $-c\sqrt{1-a}$ is sufficiently small.
\item[$(e)$] The bifurcation curve of the $2$-saddle heteroclinic loop $$HE:=\{(a,b,c)\in \Omega:b=\varphi(0,c),a=0,c<0\}.$$
\item[$(f)$] The bifurcation surface of the upper saddle connection $$SC_{1}:=\{(a,b,c)\in \Omega:b=\psi_{1}(a,c),0<a<1\},$$
where $\psi_{1}(0,c)=\varphi(0,c),\psi_{1}(a,c)>\varphi(a,c)$ for $0<a<1$ and $c<0$, and $\psi_{1}\in C^{0}$ is
increasing with respect to $a$. In addition, there exists a unique constant $a^{*}$ satisfying $\psi_{1}(a^{*},c)=-c\sqrt{1-{a^{*}}^{2}}$.
\item[$(g)$] The bifurcation surface of the lower saddle connection $$SC_{2}:=\{(a,b,c)\in \Omega:b=\psi_{2}(a,c),0<a<a_{*}<1,c<0\}$$ for a unique constant $a_{*}\in(0,1)$, where $\psi_{2}(0,c)=\varphi(0,c), \psi_{2}(a_{*},c)=0,0<\psi_{2}(a,c)<\varphi(a,c)$ for $0<a<a_{*}$ and $c<0$, and $\psi_{2}\in C^{0}$ is decreasing with respect to $a$.
\end{enumerate}
\item[(ii)] System \eqref{eq3} has exactly one limit cycle if $0\leq a<1,\varphi(a,c)<b<-c\sqrt{1-a^{2}},c<0$; and no limit cycle if one of the following conditions holds:
\begin{enumerate}
\item[$(H_{1})$] $a\geq 1$;
\item[$(H_{2})$] $0\leq a<1,b\geq -c$;
\item[$(H_{3})$] $0\leq a<1,0<b\leq \varphi(a,c),c<0$;
\item[$(H_{4})$] $0\leq a<1,-c\sqrt{1-a^{2}}\leq b<-c,c<0$.
\end{enumerate}
\end{enumerate}
\end{Theorem}

The limit cycle obtained in the above theorem is contained in the strip $-\pi-\arcsin a<x<\pi+\arcsin a$, that is, this limit cycle is homotopy to zero on the cylinder $\mathbb S^1\times \mathbb R$. In general, we call this limit cycle is contractible or first kind. From the physical background of the Josephson equation, it also has non-contractible or the second kind of limit cycles, which can be expressed as $y=y(x)$, where $y(x)$ is a smooth function with period $2\pi$. Up to now, as far as we know, the second kind of limit cycles of the Josephson equation   have not yet been considered, and studying the existence and number of such limit cycles is the main purpose of this paper.

For the second kind of limit cycle, we have
\begin{Theorem}\label{mainth2}
The Josephson equation \eqref{eq3} has at most two limit cycles of the second kind, and the upper bound is sharp.
\end{Theorem}

It is also interesting to consider the coexistence  of two  kinds of limit cycles.
\begin{Theorem}\label{mainth3}
For the Josephson equation  \eqref{eq3}, the sum of the number of limit cycles of the first and second kind is at most two, and the upper bound is sharp. Furthermore, if there are exactly two limit cycles, their configuration is $(i,j)=(0,2)$ or $(1,1)$, where $i$ (resp. $j$) represents the number of limit cycles of the first (resp. second) kind.
\end{Theorem}

To study the second kind of limit cycles, we will make the transformation $y=\frac{1}{u}$, then replace $u$ by $y$, and transform  system \eqref{eq3} to  the following Abel equation
\begin{equation}\label{eq4}
 \frac{dy}{dx}=f(x)y^2+g(x)y^3,
\end{equation}
where $f(x)= b+c\cos x, g(x)= \sin x-a, x\in [0,2\pi]$. To make the transformation $y=\frac{1}{u}$, we claim that the second kind of limit cycle does not intersect with $y=0$ in system \eqref{eq3}, and the detailed proof is given in Lemma \ref{lem2.2}. Thus we only need to consider the non-zero limit cycles  of Abel equation \eqref{eq4}. Here we call a solution $y=y(x)$ of equation \eqref{eq4} is a limit cycle if it is an isolate periodic solution with period $2\pi$.

The problem of the limit cycles of Abel equation is an important part of a big program, which is to study the limit cycle of the following non-autonomous differential equation
\begin{equation}\label{eq5}
 \frac{dx}{dt}=\displaystyle\sum_{i=0}^{n} a_{i}(t)x^{i},
\end{equation}
 where $x\in \mathbb R,t\in[0,2\pi]$, $a_{i}: \mathbb R\rightarrow \mathbb R,i=1,...,n,$ are analytic functions of period $2\pi$. There is a problem, similar to Hilbert's $16$th problem, is to ask for the maximum number of limit cycles of equation \eqref{eq5}. Until now, lots of excellent results have been obtained.

When $n\leq 2$, one can easily show that at most two limit cycles can occur. But when $n=3$, Lins Neto first pointed out in \cite{lins1980number} that the number of limit cycles of equation \eqref{eq5} is unbounded (also see \cite{panov1999diversity}). Later, Gasull and Guillamon \cite{gasull2006limit} proved that this result can be easily extended to the case of $n>3$. That is to say, to obtain an upper bound on the number of limit cycles, we need to add some conditions. The most classical condition is that certain coefficient function  or certain combination of coefficient functions has definite sign, we refer the reader to papers \cite{ alvarez2015limit,alvarez2007new,bravo2009limit,gasull2006limit,gasull1990limit,huang2017estimate,huang2020geometric,huang2012periodic,ilyashenko2000hilbert,lloyd1979note} and therein.

In particularly, when $n=3$ and $a_{0}(t)=a_{1}(t)=0$,  equation \eqref{eq5} becomes Abel equation we need to study. For this case, {\'A}lvarez, Gasull and Giacomini in \cite{alvarez2007new} gave a good criterion that if there exist real numbers $p$ and $q$ such that $pa_{2}(t)+qa_{3}(t)$ has definite sign, then equation \eqref{eq5} has at most one nonzero periodic solution, which is hyperbolic if it exists. Here, the function $h(x)$ has definite sign on the interval $J$, which means that whether $h(x)\geq 0$ or $h(x)\leq 0$ and $h(x)\not\equiv0$.

It is worth pointing out that the problem of limit cycles of equation \eqref{eq5} is extremely difficult. Even for Abel equation, we do not know whether $3$ is the maximum number of limit cycles for the following Abel equation
 \begin{equation}\label{openproblem} \frac{dx}{dt}=(a_0+a_1\sin t+a_2\cos t)x^3+(b_0+b_1\sin t+b_2\cos t)x^2. \end{equation}
In fact, Gasull \cite{gasull2021some} raised $33$ open problems and the above problem is just Problem $6$ of $33$ open problems. Since $x=0$ is always a periodic solution, Problem $6$ is to ask whether $2$ is the maximum number of non-zero limit cycles of equation \eqref{openproblem}. Notice that equation \eqref{eq4}
is a special form of equation \eqref{openproblem}, our results can be viewed as a step to attack this open problem.

Obviously, each non-zero periodic solution of equation \eqref{eq4} does not intersect with $y=0$. Therefore, we will study the number of limit cycles in regions $y>0$ and $y<0$ separately. In addition, from the analysis in \cite{chen2020global}, without loss of generality, we assume that $a\geq 0,b>0$ and $c\in \mathbb R.$

The rest of the paper is organized as follows. In Section \ref{pre}, some preliminaries which are useful tools for proving our theorems are given. In Section \ref{existence}, we discuss the existence and number of limit cycles of equation \eqref{eq4} when the parameters satisfy different conditions. The proofs of the main results and numerical example are given in the last section.

\section{Preliminaries\label{pre}}

In this section, we need to present some preliminary results that will be used in the remainder of this paper.

To study the limit cycles of equation \eqref{eq4}, the Poincar\'{e} map is a classic tools. Here we only discuss the case $y>0$, since the case $y<0$ is almost the same. Assume that $y(x, x_0, y_0)$ is the solution of equation \eqref{eq4} with initial condition $y(0, x_0, y_0)=y_0$. By the continuous dependence of solutions on the initial
value, we have $y(x, x_0, y_0)$ is well defined on $[x_0, x_0+2\pi]$ when $y_0$ is small. Denote by $I_{x_0}$ the set $\{y_0 | y(x, x_0, y_0)\mbox { is well defined on}
[x_0, x_0+2\pi]\}$ and let $y^M\triangleq \sup I_{x_0}$ ($y^M$ may be $+\infty$), then for $y_0\in (0, y^M)$, we can define the Poincar\'{e} map $G$ such that $G(y_0)=y(x_0+2 \pi, y_0)-y_0$. There is a one-one map between the isolated zero of $G$  and the positive limit cycle  of  equation \eqref{eq4}.

\subsection{Stability of the zero solution\label{Sta}}

The idea of the following procedure for the stability of the zero solution is classic, see for example \cite{lloyd1973number,lloyd1983small}.

In the Poincar\'{e} map, we choice $x_0=0$, and denote by $y(x,y_0)$ the solution of equation \eqref{eq4} with initial condition $y(0,y_0)=y_0$. Suppose that
\begin{equation}\label{eq6}
  y(x,y_0)=y_1(x)y_0+y_2(x)y_0^2+\ldots,
\end{equation}
where $y_1(0)=1,y_i(0)=0, i=2,3 , \ldots$, since  $y(0,y_0)\equiv y_0$.

Substituting \eqref{eq6} into equation \eqref{eq4}, we can easily obtain
\begin{equation}\label{eq7}
  \frac{dy_1(x)}{dx}y_0+\frac{dy_2(x)}{dx}y_0^2+...=F_2(x)y_0^2+F_3(x)y_0^3+...
\end{equation}
and
$$y_1(x)=1,$$
$$y_2(x)=\int_{0}^{x} f(t)\, dt\,,$$
$$y_3(x)=\int_{0}^{x}(2f(t)y_2(t)+g(t))\, dt\,,$$
$$...$$
where $F_2(x)={y_1^2(x)}f(x),F_3(x)=2f(x)y_1(x)y_2(x)+g(x){y_1^3(x)},...$

Substituting the expression of $f(x)$ and $g(x)$ into  equation \eqref{eq7},  we have
\begin{equation}\label{eq8}
  G(y_0) =\sum_{i=2}^{\infty} G_i(2 \pi)y_0^i,
\end{equation}
where $G_2=2\pi b, G_3=2\pi(2\pi b-a), G_4=2\pi b(\frac{4}{3}\pi ^2b^2+\pi b-4\pi a-2)+3\pi c$.

When $G_2=G_3=G_4=0$, then $a=b=c=0$ and $G(y_0)\equiv 0$, equation \eqref{eq4} has no limit cycle. Therefore, it is more important
to consider the case that $G_2, G_3$ and $G_4$ are not all zero, now the sign of the first non-zero term
will determine the stability of the zero solution, as indicated in Table \ref{Tab-zero}.
\begin{table}[!hbp]
 \caption{Stability of the zero solution of equation \eqref{eq4}}\label{Tab-zero}
  \centering
 {\footnotesize
 \setlength{\tabcolsep}{6mm}
  \begin{tabular}{|c|c|}
  \toprule
{Classification of parameters}   & Stability of the zero solution\\
\hline
  $b>0$  &  upper unstable  and lower stable\\
\hline
$b=0, a>0$  &  upper stable  and lower  stable\\
\hline
$a=b=0, c>0$  &upper unstable  and lower stable\\
\hline
$a=b=0, c<0$  &upper stable  and lower unstable\\

\bottomrule 
\end{tabular}}
\end{table}
\subsection{Stability at infinity\label{infinity}}

The Poincar\'{e} map $G(y_0)$ is defined on $(0, y^M)$, thus in this subsection we will study the stability of $G(y_0)$ near the endpoint $y^M$. Since the stability
is independent of the choice of $x_0$, we will choice $x_0$ such that $y^M=+\infty$. We will say that $y=+\infty$ is stable (resp. unstable) if $G(y_0) > ({\rm resp.}<) 0$ for sufficiently big $y_0$ (even if in general $y=+\infty$ is not an orbit).

Before proving the stability at infinity of equation \eqref{eq4}, we show that the transformation $y \rightarrow  {1}/{y}$ is reasonable.

\begin{figure}[!htbp]
\centering
\psfrag{1}{A}
\psfrag{2}{B}
\psfrag{3}{C}
\psfig{file=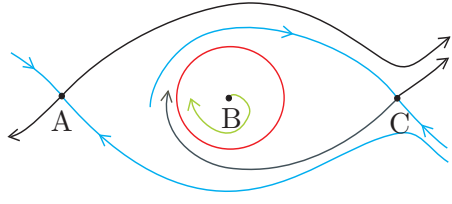,height=1.0in,width=2.3in}
\begin{center}
\caption{Phase portrait of $(a,b,c)\in S_{3}$.}
\label{Fig. 2}
\end{center}
\end{figure}

\begin{Lemma}\label{lem2.2}
The second kind of limit cycle does not intersect with $y=0$ in system \eqref{eq3}. \end{Lemma}
\noindent\bpf For system \eqref{eq3}, when $a \geq 1$, along the straight line $y=0$, $\frac{dy}{d \tau}=a-\sin x$ has the fixed sign, then any solution $y=y(x)$ which starts from $(x_0, 0)$ cannot intersect the line  $y=0$ again, thus this  solution $y=y(x)$ cannot be a period function. When $0\leq a<1$, the result will be obtained from the phase portraits obtained in  \cite{chen2020global}. From Theorem $2$ in \cite{chen2020global}, system \eqref{eq3} has three equilibria $A(-\pi-\arcsin a,0), B(\arcsin a,0)$ and $C(\pi-\arcsin a,0)$ for $0\leq a<1$,  where both $A$ and $C$ are saddles, and $B$ is an anti-saddle. Here we suppose that  the parameter $(a,b,c)\in S_{3}$, and the associated phase portrait is Fig.~\ref{Fig. 2}. If the conclusion of Lemma \ref{lem2.2} does not hold, then we can find a second kind of limit cycle $y=y(x)$, which satisfying that $y(x_M)=y(x_M+2\pi)=0$. Without loss of generality, we can assume $x_M\in  (-\pi-\arcsin a, \pi-\arcsin a)$, then this limit cycle connects the points $(x_M, 0)$ and $(x_M+2\pi, 0)$. From the phase portrait, it is easy to see that this limit cycle must intersect with the stable or unstable manifold of the saddle $C$, which is a contradiction. For other $(a,b,c)$, the proofs are almost the same, we omit them. \epf

\begin{figure}[!htbp]
\centering
\setcounter{subfigure}{0}
\subfigure[phase portrait of $(a,b,c)\in S_{1}$]
{\psfrag{1}{A}
\psfrag{2}{B}
\psfrag{3}{C}
\psfrag{4}{$l_{1}$}
\psfrag{5}{$l_{2}$}
\psfrag{6}{\textcolor{cyan}{$l_{3}$}}
\psfrag{7}{\textcolor{cyan}{$l_{4}$}}
\psfig{file=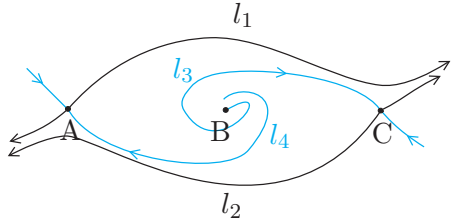,height=1.2in,width=2.3in}}\\
\centering
\subfigure[stability of $y=+\infty$ for $(a,b,c)\in S_{1}$]
{\psfrag{1}{$0$}
\psfrag{2}{$\pi/2$}
\psfrag{3}{$\pi$}
\psfrag{4}{$3\pi/2$}
\psfrag{5}{$2\pi$}
\psfrag{6}{$x$}
\psfrag{7}{$+\infty$}
\psfrag{8}{$l_{1}$}
\psfrag{9}{\textcolor{cyan}{$l_{3}$}}
\psfrag{0}{$y$}
\psfig{file=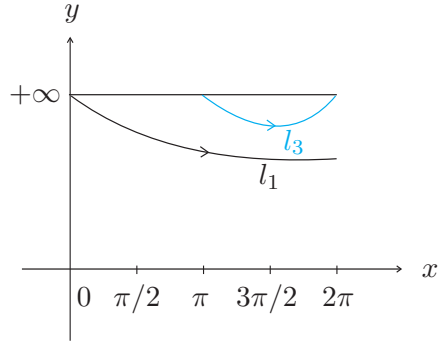,height=1.8in,width=2.3in}}
\begin{center}
\caption{Stability at infinity for $(a,b,c)\in S_{1}$.}
\label{Fig. 3}
\end{center}
\end{figure}

Then we begin to study  the stability at infinity of equation \eqref{eq4} and our proof strongly depends on  the phase diagram Fig. 2 of the paper \cite{chen2020global}. Since all the proofs are similar, we only give the detailed proof of the case $(a,b,c)\in S_{1}$, whose phase portrait is as shown in Fig.~\ref{Fig. 3}(a). Notice that the solution curve $l_1$, which is the unstable manifold of the saddle $A$, lies always in the region $y>0$. After the transformation $y\rightarrow 1/y$, $l_1$ becomes
a curve, which starts from $(-\pi-\arcsin a, +\infty)$ and is well defined for all $x>-\pi-\arcsin a$.

We set $x_0=-\pi-\arcsin a$, then for any $y_0>0$, the solution curve $y=y(x, x_0, y_0)$ is below the curve $l_1$, thus it
is well defined for all $x>x_0$. At the same time, it is easy to check that the Poincar\'{e} map
$$
G(y_0)=y(x_0+2\pi, x_0, y_0)-y_0<0,
$$
for $y_0\gg 1$, so the $y=+\infty$ is unstable, see Fig.~\ref{Fig. 3}(b).

For most cases, the stability of $\pm\infty$ can be obtained, as well as our conclusion of stability at infinity, as shown in the following table, and the detailed meanings of $a_{*},a^{*},\psi_{1}(a,c),\psi_{2}(a,c),$ $\varphi(a,c)$ can be find in \cite{chen2020global}. Furthermore, from Proposition $8$ in \cite{chen2020global}, we also have $0<\psi_{2}(a,c)<\varphi(a,c)<-c\sqrt{1-a^{2}},c<0$ for $0<a<1$ and $0<\psi_{1}(a,c)\leq-c\sqrt{1-a^{2}},c<0$ for $0\leq a\leq a^{*}$. In the table, for some cases, there is a heteroclinic orbit connecting the saddle $A$ and $C$, the above method does not work and we cannot obtain the stability of $\pm\infty$ in this way. By the way, the heteroclinic orbit connecting the saddle $A$ and $C$ can be looks as a function $y=y(x)$ with period $2\pi$ in $x$, but since this heteroclinic orbit connects two singularities,  it is not included in the limit cycle of this paper.
\begin{table}[!hbp]
\begin{center}
 \caption{Stability of system \eqref{eq4} at infinity}\label{Tab 1}
\end{center}
  \centering
 {\footnotesize
 \setlength{\tabcolsep}{4.8mm}
  \begin{tabular}{|c|c|c|c|}
  \toprule
{Region} & {Classification of parameters}   & {Stability of $y=+\infty$}  &  {Stability of $y=-\infty$}\\
\hline
  $S_{1}$ &$0<b<\psi_{2}(a,c),0<a<a_{*},c<0$  &unstable      &stable\\
\hline
$S_{2}$ &\tabincell{c}{$\psi_{2}(a,c)<b<\varphi(a,c),0<a<a_{*},c<0$\\$0<b<\varphi(a,c),a_{*}\leq a<1,c<0$}   &unstable  &unstable\\
\hline
$S_{3}$ &\tabincell{c}{$\varphi(a,c)<b<\psi_{1}(a,c),0<a<a^{*},c<0$\\$\varphi(a,c)<b<-c\sqrt{1-a^{2}},a^{*}\leq a<1,c<0$}           &unstable  &unstable\\
\hline
$S_{4}$ &$\psi_{1}(a,c)<b<-c\sqrt{1-a^{2}},0<a<a^{*},c<0$ &stable  &unstable\\
\hline
$S_{5}$ &\tabincell{c}{$b\geq-c\sqrt{1-a^{2}},0<a<a^{*},c<0$\\$b>\psi_{1}(a,c),a^{*}\leq a<1,c<0$\\$b> \max \{\psi_{1}(a,c),0\},0<a<1,c\geq0$} &stable  &unstable\\
\hline
$S_{6}$ &\tabincell{c}{$-c\sqrt{1-a^{2}}\leq b<\psi_{1}(a,c),a^{*}<a<1,c<0$\\$0<b<\psi_{1}(a,c),0<a<1,c\geq0$} &unstable  &unstable\\
\hline
$S_{7}$ &$a>1$  &unstable  &unstable\\
\hline
$HL$ &$b=\varphi(a,c),0<a<1,c<0$  &unstable  &unstable\\
\hline
$SC_{11}$ &\tabincell{c}{$b=\psi_{1}(a,c),a^{*}\leq a<1,c<0$\\$b=\psi_{1}(a,c),0<a<1,c\geq0$}  &-  &unstable\\
\hline
$SC_{12}$ &$b=\psi_{1}(a,c),0<a<a^{*},c<0$  &-  &unstable\\
\hline
$SC_{2}$ &$b=\psi_{2}(a,c),0<a<a_{*}<1,c<0$  &unstable  &-\\
\hline
$P_{1}$ &$a=1,b=\psi_{1}(1,c),c<0$  &-  &unstable\\
\hline
$SN_{1}$ &$a=1,b>\psi_{1}(1,c)$  &stable  &unstable\\
\hline
$SN_{2}$ &$a=1,0<b<\psi_{1}(1,c)$  &unstable  &unstable\\
\hline
$BT$ &$a=1,b=0,c<0$ &unstable  &unstable\\
\hline
$HLC$&$0<a<1,b=c=0$ &unstable  &unstable\\
\hline
$HE$&$a=0,b=\varphi(0,c),c<0$ &-  &-\\
\bottomrule 
\end{tabular}}
\end{table}
\subsection{The multiplicity of limit cycles\label{The}}

The Poincar\'{e} map $G(y_0)$  can be used to characterize the multiplicity of limit cycles of  Abel equation \eqref{eq4}:
\begin{enumerate}
\item[$(1)$] When $G(y_0)=0,G^{'}(y_0) \neq 0$, $y(x, x_0, y_0)$ is a hyperbolic limit cycle starting from $(x_0, y_0)$.
\item[$(2)$] When $G(y_0)=G^{'}(y_0)=...=G^{(k-1)}(y_0)=0,G^{(k)}(y_0) \neq 0$, the multiplicity of limit cycle $y(x, x_0, y_0)$ is $k$.
\end{enumerate}

To decide the multiplicity of the limit cycle, we need the following lemma, which is proved in \cite{lloyd1979note} and gives the first two derivatives of $y(2\pi, x_0, y_0)$.
\begin{Lemma}\label{lem1}
 Consider the differential equation ${\frac{dy}{dx}}=S(y,x)$, and denote by $\varphi(y_0)=y(2\pi, x_0, y_0)$, then
 \begin{enumerate}
\item[$(i)$]  $\varphi^{'}(y_0)=\exp\bigg[\displaystyle{\int_{x_0}^{x_0+2\pi} \frac{\partial S}{\partial y}(y(x, y_0),x)\, dx}\bigg];$
\item[$(ii)$] $\varphi^{''}(y_0)=\varphi^{'}(y_0)\bigg[\displaystyle{\int_{x_0}^{x_0+2\pi} \frac{\partial^{2} S}{\partial y^{2}}(y(x, y_0),x)}\cdot \exp{\displaystyle{\bigg\{\int_{x_0}^{x} \frac{\partial S}{\partial y}(y(u, y_0),u)} \, du \bigg\}}\, dx\bigg].$
\end{enumerate}
\end{Lemma}
\begin{Corollary}\label{cor1}
For equation \eqref{eq4}, if $y(x, x_0, y_0)$ is the periodic orbit starting from $(x_0, y_0)$ and $\int_{x_0}^{x_0+2 \pi}(2f(x)y(x)+3g(x)y^2(x))\,dx\,>0$ (resp. $<0$), then $y(x, x_0, y_0)$ is an unstable (resp. stable) limit cycle.
\end{Corollary}
\noindent\bpf The periodic orbit $y(x, x_0, y_0)$ shows that $G(y_0)=0$. It follows from Lemma \ref{lem1} that
 \begin{equation}\label{eq9}
  G^{'}(y_0)=\varphi^{'}(y_0)-1=\exp\left(\int_{x_0}^{x_0+2 \pi} (2f(x)y(x)+3g(x)y^2(x))\, dx\right)-1.
\end{equation}
 Therefore, if $\int_{x_0}^{x_0+2 \pi} (2f(x)y(x)+3g(x)y^2(x))\, dx>0$ (resp. $<0$), i.e. $G^{'}(y_0)>0$ (resp. $<0$), then $y(x, x_0, y_0)$ is an unstable (resp. stable) limit cycle.\epf

\subsection{Rotated vector fields\label{Rotated}}

The following results are the main tools to prove the theorem in this paper, which are obtained by the comparison theorem, since its proofs and conclusions are similar to that of the classical rotated vector fields introduced by Duff \cite{duff1953limit} and Perko \cite{perko2001differential}, so we recall the concept of rotated vector fields (it is introduced in \cite{bravo2015stability,bravo2023stability,han2018theory}) and omit the proofs.

\begin{Definition}\label{def1}
Consider the equations
\begin{equation}\label{eq10}
  \frac{dy}{dx}=f(x,y,\lambda),
\end{equation}
with parameter $\lambda\in \mathbb R$, where $x\in[0, 2\pi], y\in \mathbb R$ and $f(x, y, \lambda)$ is $2\pi$-periodic function with respect to $x$, which is called the rotated vector fields if  $\frac{\partial{f(x,y,\lambda)}}{\partial{\lambda}}>0$.
\end{Definition}
\begin{Proposition}\label{pro1}
Suppose that equation \eqref{eq10} be rotated vector fields for $\lambda\in \mathbb R$, the limit cycles of $\eqref{eq10}_{\lambda_1}$ and $\eqref{eq10}_{\lambda_2}$ cannot intersect if  $\lambda_1\neq \lambda_2$.
\end{Proposition}

\begin{Proposition}\label{pro2}
The stable (resp. unstable) and hyperbolic limit cycle of equation \eqref{eq10} will expands (resp. contracts ) as ${\lambda}$ increase. Here we call that a limit cycle
$y=y(x, \lambda)$ expands (resp. contracts ) if $y(x, \lambda)$ is bigger (resp. smaller) for each $x\in \mathbb R$.
\end{Proposition}

\begin{Proposition}\label{pro3} For semi-stable limit cycle of multiplicity two, when ${\lambda}$ increases, the upper unstable limit cycle will disappear, and the upper stable limit cycle will split into two hyperbolic limit cycles .
\end{Proposition}

\section{The existence and the number of limit cycles\label{existence}}

In this Section, we focus on the existence and the number of limit cycles of equation \eqref{eq4}.

It follows from Theorem A of paper \cite{alvarez2007new} that when ${\xi}g(x)+f(x)$  has definite sign in $[0,2\pi]$ for some ${\xi}\in \mathbb R$, equation \eqref{eq4} has at most one limit cycle. Moreover, when the limit cycle exists, it is hyperbolic. Thus we need to consider the function:
          $${\xi}g(x)+f(x)={\xi}(\sin x-a)+b+c \cos x=\sqrt{{\xi}^2+c^2}\sin(x+\theta)+b-{\xi}a.$$
When it has definite sign in $[0,2\pi]$ i.e. for some $\xi$, $\xi^2+c^2\leq (b-\xi a)^2$, which is equivalent to one of the following three conditions hold:
\begin{enumerate}
\item[$(A_{1})$] $a\geq1,b>0$;
\item[$(A_{2})$] $0\leq a<1,b\geq\lvert c\rvert$;
\item[$(A_{3})$] $0\leq a<1,\lvert c\sqrt{1-a^2}\rvert\leq b<\lvert c\rvert$.
\end{enumerate}
Recall that we have assumed that $a\geq 0,b>0$ and $c\in \mathbb R.$

And when ${\xi}g(x)+f(x)$ has no definite sign in $[0,2\pi]$, i.e. $0<a<1,0<b<\lvert c\sqrt{1-a^2}\rvert$, the conclusion of Theorem A in Paper \cite{alvarez2007new} does not work, and more limit cycles can appear.

\subsection{The case: \texorpdfstring{${\xi}g(x)+f(x)$}{}  has definite sign}
 For the case of ${\xi}g(x)+f(x)$  has definite sign in $[0,2\pi]$, the key point to prove the existence of limit cycles is to discuss the stability of $y=0$ and $y=\infty$. Concretely, we have the following lemma.
\begin{Lemma}\label{lem3.1}
For equation \eqref{eq4}, if the stabilities of $y=0$ and $y=\infty$ are same (resp. opposite), then the number of limit cycles is odd (resp. even) (taking into account multiplicities).
\end{Lemma}
\bpf We first prove the case outside of the brackets and $y>0$. Assume that both $y=0$ and $y=+\infty$ are unstable, then the Poincar\'{e} map $G(s_{1})>0$ and $G(s_{2})<0$
where $0<s_1 \ll 1 \ll s_2$. One could obtain that  the number of zeros of $G(y)=0$ on $(s_{1},s_{2})$ must be odd (taking into account multiplicities), that is, the number of limit cycles of equation \eqref{eq4} is odd (taking into account multiplicities). The proofs of other cases are similar, and we omit them.\epf

Specially, when ${\xi}g(x)+f(x)$  has definite sign in $[0,2\pi]$, if the stabilities of $y=0$ and $y=\infty$ are same (resp. opposite), then  equation \eqref{eq4} have $1$ (resp. $0$) limit cycle.

\subsubsection {Limit cycles in the region  \texorpdfstring{$y>0$}{}\label{l}}

\begin{Theorem}\label{th3.1}
Equation \eqref{eq4} has exactly one limit cycle in $y>0$ if $a>1,b>0$ or $a=1,0<b<\psi_{1}(1,c)$; and no limit cycle in $y>0$ if $a=1,b\geq \psi_{1}(1,c)$.
\end{Theorem}
\noindent\bpf It can be seen from Table \ref{Tab 1} that $y=+\infty$ is unstable if $a>1$ or $a=1,0<b<\psi_{1}(1,c)$, hence we get that there is exactly one limit cycle for $a>1,b>0$ or $a=1,0<b<\psi_{1}(1,c)$ in $y>0$. In addition, $y=+\infty$ is stable if $a=1,b>\psi_{1}(1,c)$, so similarly, one could obtain that there is no limit cycle for $a=1,b>\psi_{1}(1,c)$ in $y>0$.

When $a=1,b=\psi_{1}(1,c)$, we do not know the stability of $y=+\infty$, thus we need extra proof. We claim that there is no limit cycle. Else equation \eqref{eq4} has exactly one limit cycle for $a=1,b=\psi_{1}(1,c)$, which is hyperbolic. Let $b$ change from $\psi_{1}(1,c)$ to $\psi_{1}(1,c)+\epsilon$, where $0<\epsilon\ll1$, then this hyperbolic limit cycle still exist, that is, we find a limit cycle for $a=1,b=\psi_{1}(1,c)+\epsilon$. In fact, equation \eqref{eq4} has no limit cycle for $a=1,b>\psi_{1}(1,c)$, this is a contradiction. Thus, the claim holds, which finishes the proof.\epf

In order to determine whether the limit cycle exists for $0\leq a<1$ and $b>0$, we first analyze the relationship between the functions $b=\psi_{1}(a,c),b=\lvert c\sqrt{1-a^2}\rvert$ and $b=\lvert c\rvert$ by the bifurcation diagram Fig. $1$ of \cite{chen2020global}.

When $c>0$, from the fact that the functions $b=\psi_{1}(a,c)$ and $b=c\sqrt{1-a^2}$ are monotonically increasing and monotonically decreasing with respect to $a$, respectively, we have $\psi_{1}(a,c)< \psi_{1}(1,c)$, so the following two cases need to be considered:

Case $1$: $0<c<\psi_{1}(1,c)$

\begin{figure}[!htbp]
\centering
\setcounter{subfigure}{0}
\subfigure[$0<c<\psi_{1}(1,c)$]
{\psfrag{1}{$0$}
\psfrag{2}{$a_{*}$}
\psfrag{3}{$\bar{a}$}
\psfrag{4}{$\tilde{a}$}
\psfrag{5}{$1$}
\psfrag{6}{$a$}
\psfrag{7}{$c$}
\psfrag{8}{$1/c$}
\psfrag{9}{$b$}
\psfrag{a}{$b=c\sqrt{1-a^{2}}$}
\psfrag{b}{$b=c$}
\psfrag{c}{$b=\psi_{1}(a,c)$}
\psfig{file=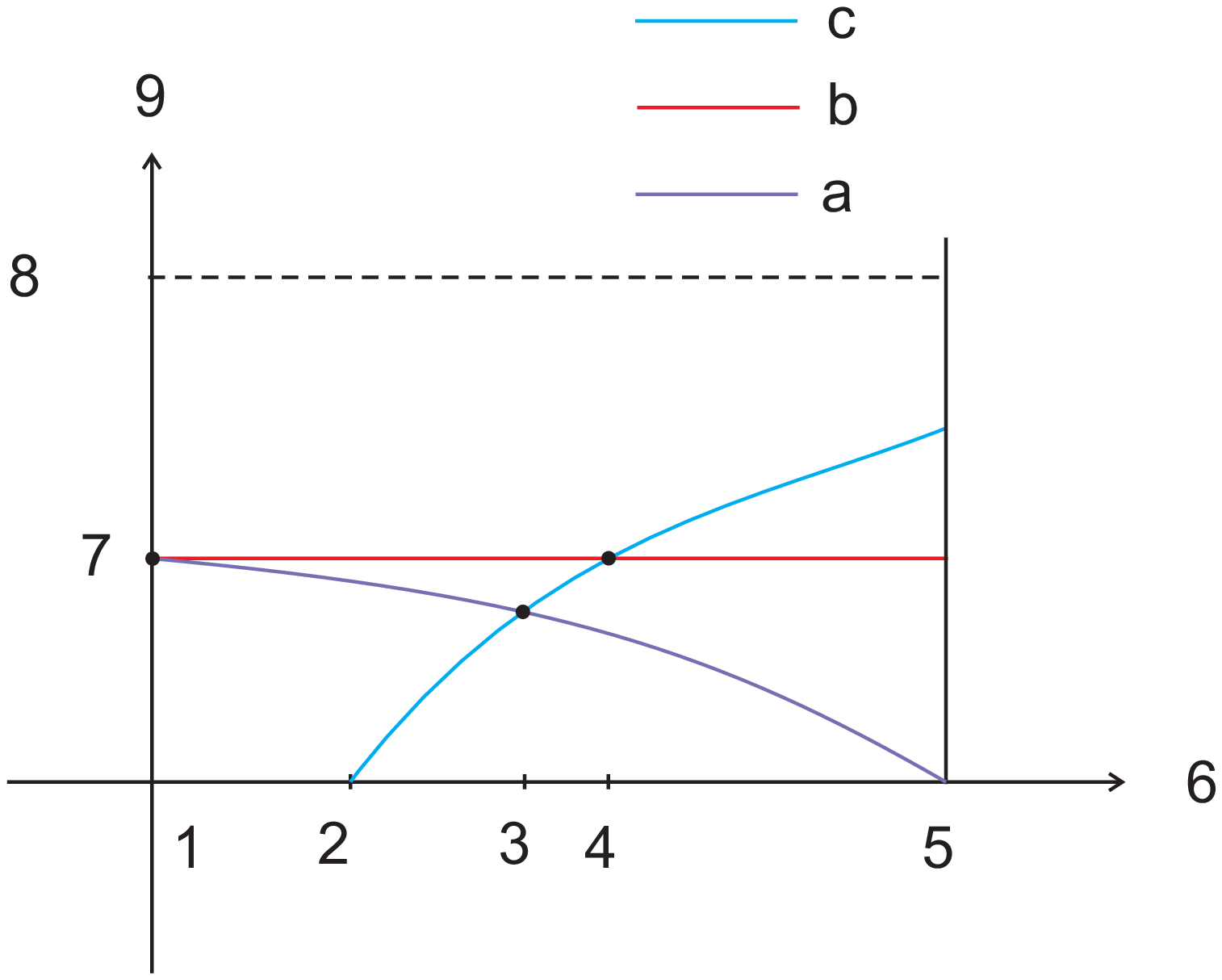,height=1.9in,width=2.4in}}
\hspace{21mm}
\centering
\subfigure[$c\geq \psi_{1}(1,c)$]
{\psfrag{1}{$0$}
\psfrag{2}{$a_{*}$}
\psfrag{3}{$\bar{a}$}
\psfrag{4}{$1$}
\psfrag{5}{$a$}
\psfrag{6}{$c$}
\psfrag{7}{$1/c$}
\psfrag{8}{$b$}
\psfrag{9}{$b=c\sqrt{1-a^{2}}$}
\psfrag{a}{$b=c$}
\psfrag{b}{$b=\psi_{1}(a,c)$}
\psfig{file=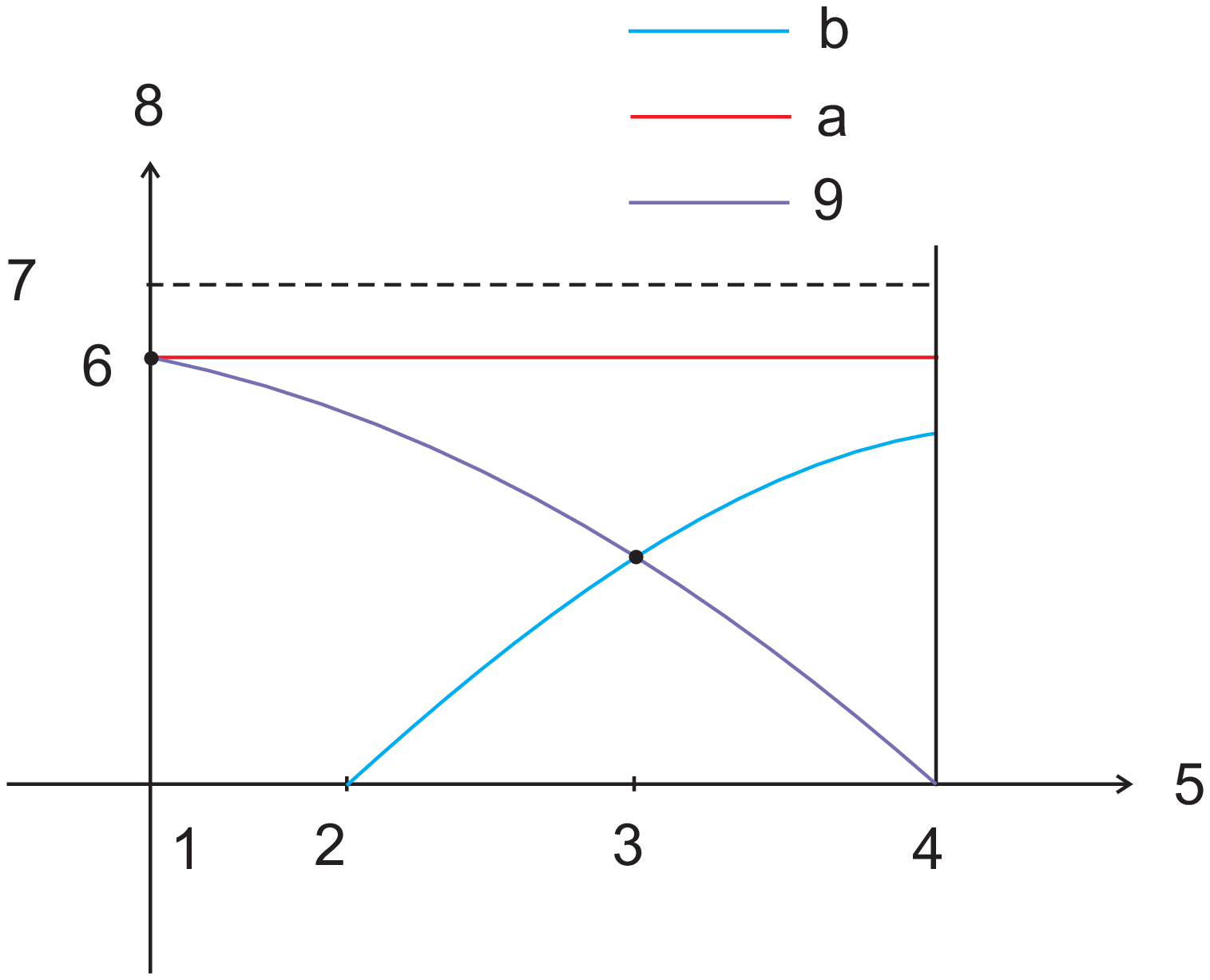,height=1.9in,width=2.4in}}
\\
\centering
\subfigure[$0<-c<\psi_{1}(1,c)$]
{\psfrag{1}{$0$}
\psfrag{2}{$a^{*}$}
\psfrag{3}{$\hat{a}$}
\psfrag{4}{$1$}
\psfrag{5}{$a$}
\psfrag{6}{$-c$}
\psfrag{7}{$b$}
\psfrag{8}{$b=-c\sqrt{1-a^{2}}$}
\psfrag{9}{$b=-c$}
\psfrag{0}{$b=\psi_{1}(a,c)$}
\psfig{file=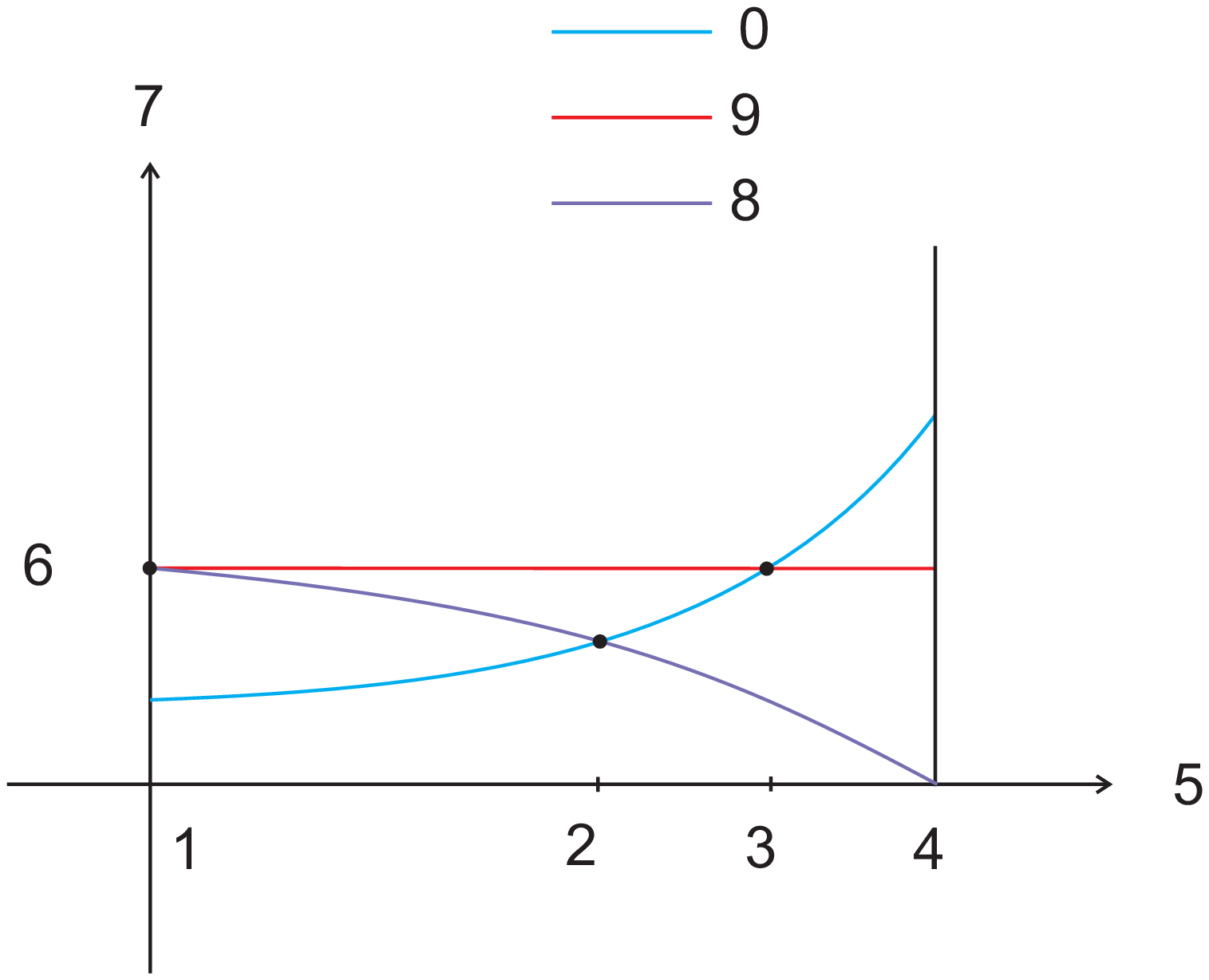,height=1.9in,width=2.4in}}
\hspace{21mm}
\centering
\subfigure[$-c\geq \psi_{1}(1,c)$]
{\psfrag{1}{$0$}
\psfrag{2}{$a^{*}$}
\psfrag{3}{$1$}
\psfrag{4}{$a$}
\psfrag{5}{$-c$}
\psfrag{6}{$b$}
\psfrag{7}{$b=-c\sqrt{1-a^{2}}$}
\psfrag{8}{$b=-c$}
\psfrag{9}{$b=\psi_{1}(a,c)$}
\psfig{file=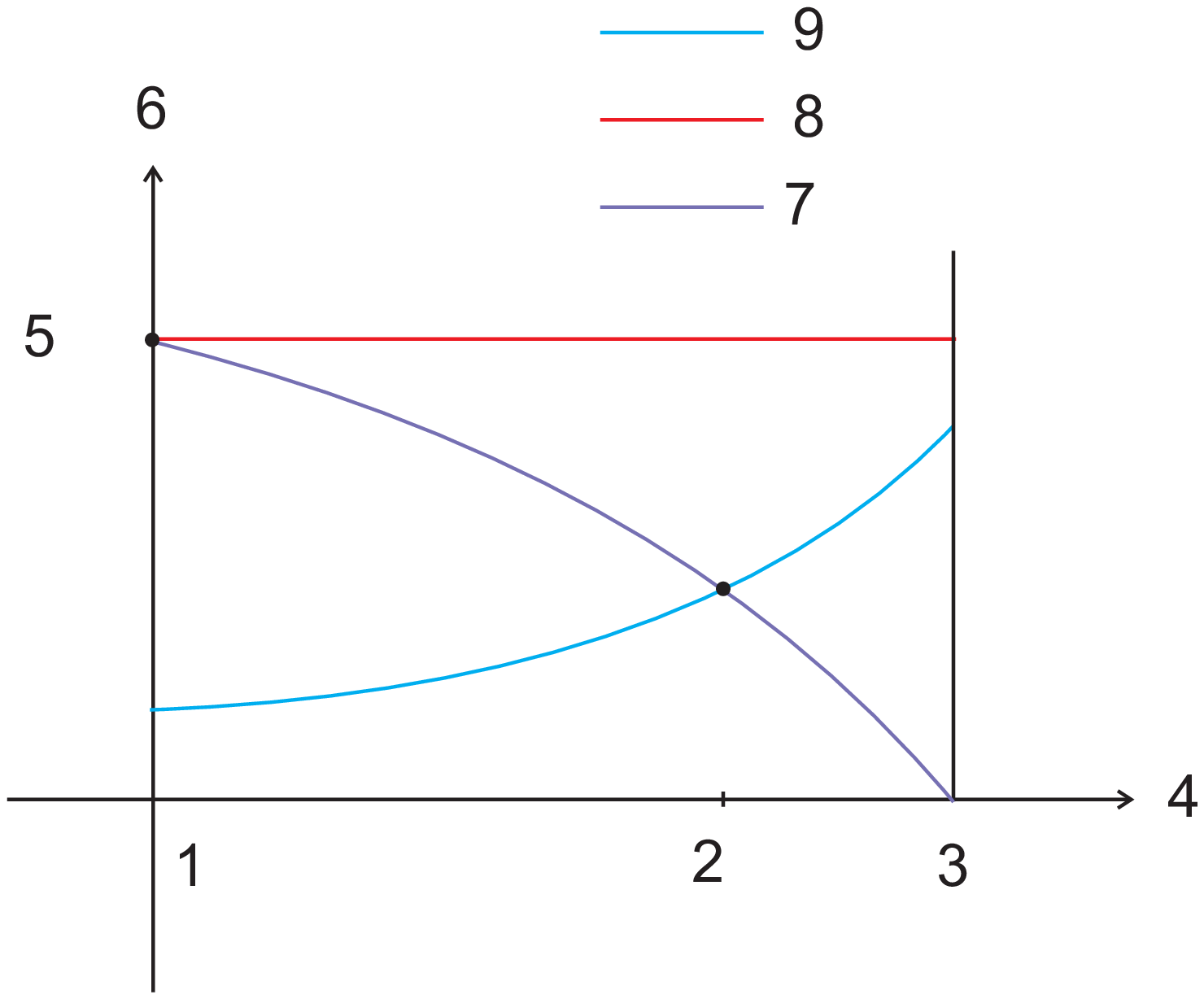,height=1.9in,width=2.4in}}\\
\centering
\subfigure[$c=0$]
{\psfrag{1}{$0$}
\psfrag{2}{$1$}
\psfrag{3}{$a$}
\psfrag{4}{$b$}
\psfrag{5}{$b=\psi_{1}(a,c)$}
\psfig{file=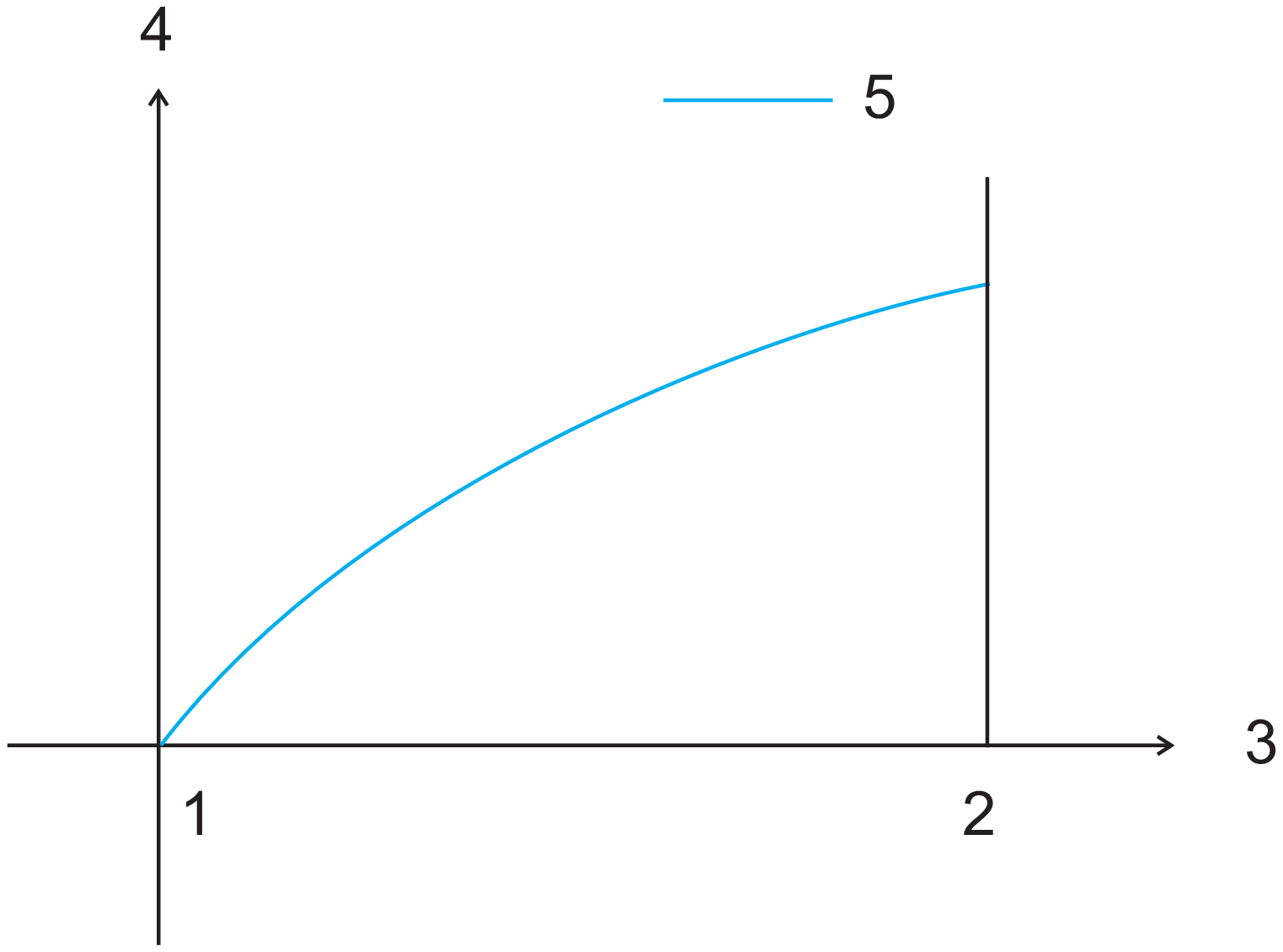,height=1.9in,width=2.4in}}
\begin{center}
\caption{Relationship between the functions $b=\psi_{1}(a,c),b=\lvert c\sqrt{1-a^2}\rvert$ and $b=\lvert c\rvert$.}
\label{Fig. 4}
\end{center}
\end{figure}

According to Fig. $1(c)$ of \cite{chen2020global}, there exists a unique $\bar{a}\in(a_{*},1)$ such that $\psi_{1}(\bar{a},c)=c\sqrt{1-\bar{a}^2}$ and a unique $\tilde{a}\in(\bar{a},1)$ such that $\psi_{1}(\tilde{a},c)=c$ (see Fig.~\ref{Fig. 4}(a)). Then we have three subcases: ${\max\{\psi_{1}(a,c),0\}}< c\sqrt{1-a^2}\leq c$ for $0\leq a< \bar{a}$; $c\sqrt{1-a^2}\leq \psi_{1}(a,c)< c$ for $\bar{a}\leq a< \tilde{a}$ and $c\sqrt{1-a^2}<c\leq \psi_{1}(a,c)$ for $\tilde{a}\leq a<1$.

Case $2$: $c\geq \psi_{1}(1,c)$

Clearly if $c\geq \psi_{1}(1,c)$, then $c>\psi_{1}(a,c)$ and $c\geq c\sqrt{1-a^2}$, so there must be a unique $\bar{a}\in(a_{*},1)$ satisfies $\psi_{1}(\bar{a},c)=c\sqrt{1-\bar{a}^2}$ (see Fig.~\ref{Fig. 4}(b)). Then we have two subcases: ${\max\{\psi_{1}(a,c),0\}}< c\sqrt{1-a^2}\leq c$ for $0\leq a< \bar{a}$ and $c\sqrt{1-a^2}\leq \psi_{1}(a,c)< c$ for $\bar{a}\leq a<1$.

When $c<0$, due to $b=\psi_{1}(a,c)$ is monotonically increasing with respect to $a$ and $b=-c\sqrt{1-a^2}$ is monotonically decreasing with respect to $a$, one could obtain that $\psi_{1}(a,c)<\psi_{1}(1,c)$, next we will discuss it in two cases:

Case $3$: $0<-c<\psi_{1}(1,c)$

It follows from Fig. $1(a)$ of \cite{chen2020global} that there is a unique $\hat{a}\in (a^{*},1)$ such that $\psi_{1}(\hat{a},c)=-c$ (see Fig.~\ref{Fig. 4}(c)), then we have three subcases:  $\psi_{1}(a,c)<-c\sqrt{1-a^2}\leq -c$ for $0\leq a<a^{*}$; $-c\sqrt{1-a^2}\leq \psi_{1}(a,c)<-c$ for $a^{*}\leq a< \hat{a}$ and $-c\sqrt{1-a^2}<-c\leq \psi_{1}(a,c)$ for $\hat{a}\leq a<1$.

Case $4$: $-c\geq \psi_{1}(1,c)$

Note that $-c\sqrt{1-a^2}\leq-c$ and $\psi_{1}(a,c)<-c$ for $-c\geq \psi_{1}(1,c)$,  then we have two subcases:  $\psi_{1}(a,c)<-c\sqrt{1-a^2}\leq -c$ for $0\leq a<a^{*}$ and $-c\sqrt{1-a^2}\leq \psi_{1}(a,c)< -c$ for $a^{*}\leq a< 1$, see Fig.~\ref{Fig. 4}(d).

When $c=0$, based on the monotonicity of the functions $b=\psi_{1}(a,c)$, we know that $\psi_{1}(a,c)\geq 0=c\sqrt{1-a^2}=c$, see Fig.~\ref{Fig. 4}(e).

\begin{Theorem}\label{th3.2}
 Suppose that $0<a<1, b>|c|$. Equation \eqref{eq4} has exactly one limit cycle in $y>0$ if and only if one of the following conditions holds:
 \begin{enumerate}
\item[$(B_{1})$] $0<a<1,0<b<\psi_{1}(a,c),c=0$;
\item[$(B_{2})$] $\tilde{a}<a<1,c\leq b<\psi_{1}(a,c),0<c<\psi_{1}(1,c)$;
\item[$(B_{3})$] $\hat{a}<a<1,-c\leq b<\psi_{1}(a,c),0<-c<\psi_{1}(1,c)$.
\end{enumerate}
And equation \eqref{eq4} has no limit cycle in $y>0$ if $0<a<1,b\geq max\{\psi_{1}(a,c),\lvert c\rvert\}$.
\end{Theorem}
\noindent\bpf From the stabilities of $y=0$ and $y=+\infty$, we can find that when one of the conditions (${B_{1}}$) to (${B_{3}}$) holds, there is exactly one limit cycle in $y>0$. At the same time, there is no limit cycle when one of the following conditions holds:
\begin{enumerate}
\item[$(B_{4})$] $0<a<1,b>\psi_{1}(a,c),c=0$;
\item[$(B_{5})$] $0<a< \tilde{a},b\geq c,0<c<\psi_{1}(1,c)$;
\item[$(B_{6})$] $\tilde{a}\leq a<1,b>\psi_{1}(a,c),0<c<\psi_{1}(1,c)$;
\item[$(B_{7})$] $0<a<1,b\geq c,c\geq \psi_{1}(1,c)$;
\item[$(B_{8})$] $0<a<1,b\geq -c,-c\geq \psi_{1}(1,c)$;
\item[$(B_{9})$] $0<a<\hat{a},b\geq -c,0<-c<\psi_{1}(1,c)$;
\item[$(B_{10})$] $\hat{a}\leq a<1,b>\psi_{1}(a,c),0<-c<\psi_{1}(1,c)$.
\end{enumerate}

For the following cases:
\begin{enumerate}
\item[$(B_{11})$] $0<a<1,b=\psi_{1}(a,c),c=0$;
\item[$(B_{12})$] $\tilde{a}\leq a<1,b=\psi_{1}(a,c),0<c<\psi_{1}(1,c)$;
\item[$(B_{13})$] $\hat{a}\leq a<1,b=\psi_{1}(a,c),0<-c<\psi_{1}(1,c)$.
\end{enumerate}
We do not know the stability of $y=+\infty$. Then by using the same idea in the proof of Theorem \ref{th3.1}, it is not difficult to derive that there is no limit cycle.

Consequently, equation \eqref{eq4} has no limit cycle for $0<a<1, b>\max\{\psi_{1}(a,c), \lvert c\rvert\}$, which implies that the conclusion of this theorem holds.\epf
\begin{Theorem}\label{th3.3} Suppose that $0<a<1, |c\sqrt{1-a^2}|\leq b<|c|$. Equation \eqref{eq4} has exactly one limit cycle in $y>0$ if $\bar{a}<a<1, c\sqrt{1-a^2}\leq b< min\{\psi_{1}(a,c),c\},c>0$ or $a^{*}<a<1, -c\sqrt{1-a^2}\leq b<min\{\psi_{1}(a,c),-c\},c<0$ holds; and no limit cycle in $y>0$ if and only if one of the following conditions holds:
\begin{enumerate}
\item[$(C_{1})$] $0<a<\bar{a},c\sqrt{1-a^2}\leq b<c,c>0$;
\item[$(C_{2})$] $\bar{a}\leq a< \tilde{a},\psi_{1}(a,c)\leq b<c,0<c<\psi_{1}(1,c)$;
\item[$(C_{3})$] $\bar{a}\leq a<1, \psi_{1}(a,c)\leq b<c,c\geq \psi_{1}(1,c)$;
\item[$(C_{4})$] $0<a<a^{*},-c\sqrt{1-a^2}\leq b<-c,c<0$;
\item[$(C_{5})$] $a^{*}\leq a<\hat{a},\psi_{1}(a,c)\leq b<-c,0<-c<\psi_{1}(1,c)$;
\item[$(C_{6})$] $a^{*}\leq a<1,\psi_{1}(a,c)\leq b<-c,-c\geq \psi_{1}(1,c)$.
\end{enumerate}
\end{Theorem}
\noindent\bpf For most cases, the conclusions can be drawn directly from the stability of $y=0$ and $y=+\infty$.
And for the following cases:
\begin{enumerate}
\item[$(C_{7})$] $\bar{a}\leq a< \tilde{a},b=\psi_{1}(a,c),0<c<\psi_{1}(1,c)$;
\item[$(C_{8})$] $\bar{a}\leq a<1,b=\psi_{1}(a,c),c\geq \psi_{1}(1,c)$;
\item[$(C_{9})$] $a^{*}\leq a<\hat{a},b=\psi_{1}(a,c),0<-c<\psi_{1}(1,c)$;
\item[$(C_{10})$] $a^{*}\leq a<1,b=\psi_{1}(a,c),-c\geq \psi_{1}(1,c)$.
\end{enumerate}
We do not know the stability of $y=+\infty$, too. Therefore, we repeat the same idea in the proof of Theorem \ref{th3.1} and obtain that there is no limit cycle. This
completes the proof of Theorem \ref{th3.3}.\epf

\subsubsection {Limit cycles in the region \texorpdfstring{$y<0$}{}\label{m}}
\begin{Theorem}\label{th3.4} Equation \eqref{eq4} has no limit cycle in $y<0$ if one of the conditions (${A_{1}}$) to (${A_{3}}$) holds. \end{Theorem}
\noindent\bpf It is easy to verified that any of the conditions (${A_{1}}$) to (${A_{3}}$) hold, equation \eqref{eq4} has at most one limit cycle and $y=0$ is stable, $y=-\infty$ is unstable. Thus,  equation \eqref{eq4} has no limit cycle in $y<0$ and the conclusion of the theorem is proven. \epf

\subsection{Some auxiliary results \label{n}}

From now on, we will consider the remainder case: $0<a<1,0<b<|c\sqrt{1-a^2}|$. The problem becomes much more difficult, new ideas must be introduced. In this subsection, we give some auxiliary results that will be used in the rest of the paper. In order to prove our theorems in the future, we will deal with a little more general case: $0\leq a<1,0\leq b<|c\sqrt{1-a^2}|$.

It can be seen from $\frac{dy}{dx}=0$ that $y(x)=-\frac{f(x)}{g(x)}$  is the isoclinic of equation \eqref{eq4}. For the sake of discussion, we will write it as $\omega(x)=-\frac{f(x)}{g(x)}$.

Below we first consider the case of $0\leq a<1,0\leq b<c\sqrt{1-a^2},c>0$ and $y>0$  and get the following results.
\begin{Lemma}\label{lem3.2}
The function $\omega(x)=-\frac{f(x)}{g(x)}$ is monotonically increasing with respect to $x$. \end{Lemma}
\noindent\bpf By simple calculations, we have
\begin{equation}\label{eq12}
\begin{split}
 \omega^{'}(x) & =(-\frac{f(x)}{g(x)})^{'}\\
   & =\frac{-f^{'}(x)g(x)+f(x)g^{'}(x)}{g^{2}(x)}\\
   & =\frac{c-ac\sin x+b\cos x}{(\sin x-a)^{2}}\\
   & =\frac{c+\sqrt{a^{2}c^{2}+b^{2}}\sin(x+\theta)}{(\sin x-a)^{2}},
    \end{split}
\end{equation}
where $\tan \theta=-\frac{b}{ac}$. If $0\leq a<1,0\leq b<c\sqrt{1-a^2},c>0$, then $c> \sqrt{a^{2}c^{2}+b^{2}}$, we can conclude that $\omega^{'}(x)>0$, i.e. the function $\omega(x)=-\frac{f(x)}{g(x)}$ is monotonically increasing with respect to $x$. \epf
\begin{Lemma}\label{lem3.3} If equation \eqref{eq4} has one limit cycle, then this limit cycle has exactly one maximum point and one minimum point in any interval of length $2\pi$. Here we look the two endpoints of the interval of length $2\pi$ as one point.  \end{Lemma}
\noindent\bpf It is enough to analyze the sign of $\frac{dy}{dx}$ in $[\frac{\pi}{2},\frac{5\pi}{2}]$.

Since
 \begin{equation}\label{eq13}
 \frac{dy}{dx}=f(x)y^2+g(x)y^3=y^2(\sin x-a)\big(y+\frac{f(x)}{g(x)}\big),
\end{equation}
 the sign of $\frac{dy}{dx}$ is determined by the signs of $\sin x-a$ and $y+\frac{f(x)}{g(x)}$.

Below we analyze their signs in $[\frac{\pi}{2},\frac{5\pi}{2}]$.

Letting $x_{1}=\arcsin a$ and $x_{2}=\pi-\arcsin a$, then $\sin x-a>0$ for $x\in (x_{1},x_{2})$ and $\sin x-a<0$ for $x\in [0,x_{1})\cup(x_{2},2\pi]$.

\begin{figure}[!htbp]
\begin{center}
\psfrag{1}{$0$}
\psfrag{2}{$x_{1}$}
 \psfrag{3}{$\pi/2$}
 \psfrag{4}{$x_{2}$}
  \psfrag{5}{$\pi$}
  \psfrag{6}{$3\pi/2$}
   \psfrag{7}{$2\pi$}
   \psfrag{8}{$5\pi/2$}
   \psfrag{9}{$3\pi$}
   \psfrag{0}{$x$}
 \psfrag{a}{$y$}
  \psfrag{b}{$\omega(x)=-f(x)/g(x)$}
   \psfrag{c}{$\sin x=a$}
   \psfrag{d}{$y(x)$}
   \psfrag{e}{}
\psfig{file=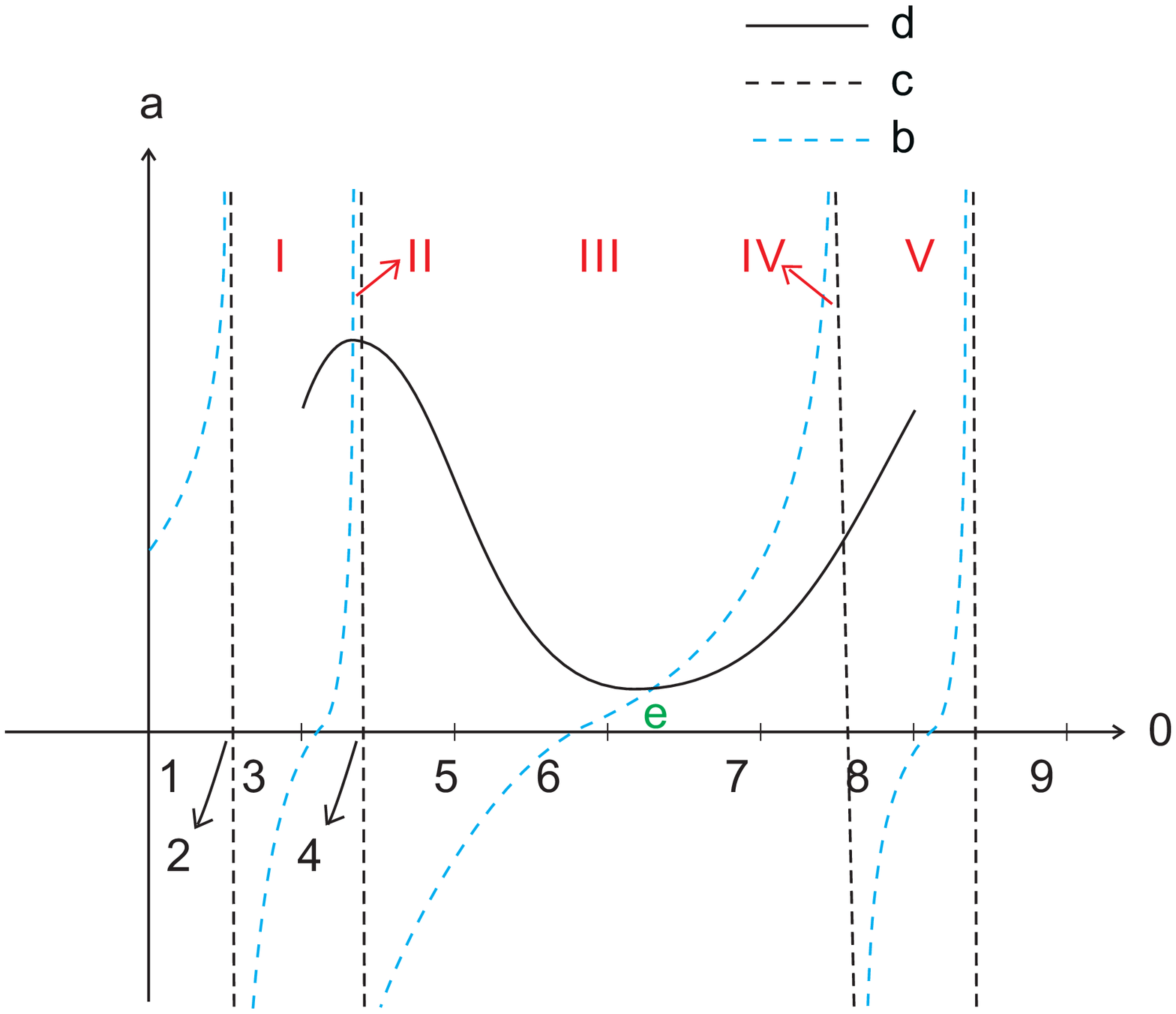,height=3.0in,width=3.2in}
\begin{center}
\caption{Solution curve of equation \eqref{eq4} when $x\in[\frac{\pi}{2},\frac{5\pi}{2}]$.}
\label{Fig. 5}
\end{center}
\end{center}
\end{figure}

Furthermore, we have known that $\omega(x)=-\frac{f(x)}{g(x)}$ is monotonically increasing and tends to $\infty$ as $x$ tends to $x_1$ and $x_2$. Therefore, the function $\omega(x)=-\frac{f(x)}{g(x)}$ divides the region $x_1<x<x_2$ into two parts $I$ (left) and $II$ (right), in which $y+\frac{f(x)}{g(x)}>0$ in part $I$ and $y+\frac{f(x)}{g(x)}<0$ in part $II$. Similarly, we can analyze the signs for parts $III,IV$ and $V$, as shown in Table \ref{Tab 2 }.
From the continuity of $y(x)$, its diagram can be plotted, as depicted in Fig.~\ref{Fig. 5}. Consequently, it is not difficult to draw the conclusion of Lemma \ref{lem3.3}. \epf
\begin{table}[!hbp]
\caption{Sign of $\frac{dy}{dx}$ in different parts}\label{Tab 2 }
  \centering
 {\footnotesize
  \setlength{\tabcolsep}{6mm}
  \begin{tabular}{|c|c|c|c|c|c|}
  \toprule
{Part} & {I}   & {II}  &  {III}& {IV}  &  {V}\\
\hline
  $\sin x-a$ &+  &+   &-  &-  &+\\
\hline
$y+\frac{f(x)}{g(x)}$ &+   &-  &+  &-  &+\\
\hline
$\frac{dy}{dx}$ &+  &-  &- &+  &+\\
\bottomrule 
\end{tabular}}
\end{table}

 \begin{figure}[!htbp]
\begin{center}
\psfrag{1}{$0$}
\psfrag{2}{$\pi/2$}
 \psfrag{3}{$x^{*}$}
 \psfrag{4}{$\pi$}
  \psfrag{5}{$3\pi/2$}
  \psfrag{6}{$x_{*}$}
   \psfrag{7}{$2\pi$}
   \psfrag{8}{$5\pi/2$}
   \psfrag{9}{$x^{*}+2\pi$}
   \psfrag{0}{$3\pi$}
 \psfrag{a}{$x$}
  \psfrag{b}{$y_{*}$}
   \psfrag{c}{$y^{*}$}
 \psfrag{d}{$y$}
  \psfrag{e}{$\omega(x)=-f(x)/g(x)$}
   \psfrag{f}{$\sin x=a$}
   \psfrag{g}{$y(x)$}
\psfig{file=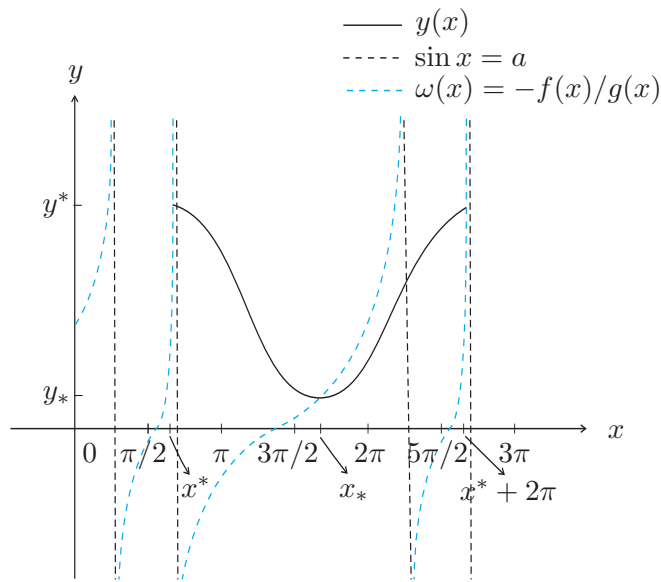,height=3.0in,width=3.3in}
\begin{center}
\caption{Solution curve of equation \eqref{eq4} when $x\in[x^{*},x^{*}+2\pi]$.}
\label{Fig. 6}
\end{center}
\end{center}
\end{figure}

\begin{Lemma}\label{lem3.4} The multiplicity of the limit cycle of equation \eqref{eq4} is at most two. Moreover, the limit cycle must be lower stable if it exists.\end{Lemma}

\noindent\bpf Assume that $y=y(x, x^*, y^*)$ is a positive limit cycle of  equation \eqref{eq4}, where $x^*$ is the unique  maximum point and  $y^*=y(x^*)$. At the same time, denote the unique minimum point by $x_{*}$ and let $y_*=y(x_*)$, as shown in Fig.~\ref{Fig. 6}.

Write any positive solution, which satisfies $y(x^{*})=y_0$, of equation \eqref{eq4} as $y(x^{*}+2\pi, x^{*}, y_0)$, and define the Poincar\'{e} map $$G(y_0)=y(x^{*}+2\pi, x^{*}, y_0)-y_0.$$

By Lemma \ref{lem1}, we have
\begin{equation}\label{eq16A}
G'(y_0)=\exp\big(\int_{x^*}^{x^*+2 \pi} (2f(x)y(x)+3g(x)y^2(x))\, dx\big)-1,
\end{equation}
and
\begin{equation}\label{eq16}
G''(y_0)=(1+G'(y_0))\int_{x^{*}}^{2 \pi+x^{*}}(2f(x)+6g(x)y(x))\exp\big(\int_{x^{*}}^{x} (2f(u)y(u)+3g(u)y^2(u))\, du\big)\, dx.
\end{equation}

If the multiplicity of the limit cycle $y=y(x, x^*, y^*)$ is bigger than $1$, then $G'(y^*)=0$. For convenience of notations, in the following proof of Lemma \ref{lem3.4}, we use
$y(x)$ instead of $y(x, x^*, y^*)$. Then we have
$$\int_{x^{*}}^{2 \pi+x^{*}}(2f(x)y(x)+3g(x)y^2(x))\, dx=0,$$
$$G''(y^*)=4I_1+2I_2,$$
where
$$I_1=\int_{x^{*}}^{2 \pi+x^{*}}g(x)y(x)\exp\big(\int_{x^{*}}^{x} (2f(u)y(u)+3g(u)y^2(u))\, du\big)\, dx, $$
 $$I_2=\int_{x^{*}}^{2 \pi+x^{*}}(f(x)+g(x)y(x))\exp\big(\int_{x^{*}}^{x} (2f(u)y(u)+3g(u)y^2(u))\, du\big)\, dx.$$

Now $y(x)\neq 0$, we can rewrite  equation \eqref{eq4} as
\begin{equation}\label{eq14}
  \frac{dy(x)}{y(x)}=\big(f(x)y(x)+g(x)y^{2}(x)\big)dx.
\end{equation}
Integrating both sides of equation \eqref{eq14} gives
 \begin{equation}\label{eq15}
 \int_{x^{*}}^{2 \pi+x^{*}} \big(f(x)y(x)+g(x)y^{2}(x)\big)\,dx=0,
\end{equation}
and
\begin{equation}\label{eq17}
  \int_{x^{*}}^{x} (f(u)y(u)+g(u)y^2(u))\, du=\ln \frac{y(x)}{y^*}.
\end{equation}

Thus $\int_{x^{*}}^{2 \pi+x^{*}}g(x)y^{2}(x)\, dx=0$.
We define a function $k(x)=\int_{x^{*}}^{x} g(u)y^2(u)\, du$. Since $k(x^{*})=k(2 \pi+x^{*})=0$, we have
\begin{equation}\label{eq18}
\begin{split}
I_1 & =\int_{x^{*}}^{2 \pi+x^{*}} g(x)y(x)\exp\big(2 \ln \frac{y(x)}{y(x^{*})}+k(x)\big)\, dx\\
   & =\frac{1}{(y^*)^2}\int_{x^{*}}^{2 \pi+x^{*}} g(x)y^{3}(x)\exp k(x)\, dx\\
   & =\frac{1}{(y^*)^2}\int_{x^{*}}^{2 \pi+x^{*}} y(x) d \exp k(x)\\
   & =-\frac{1}{(y^*)^2}\int_{x^{*}}^{2 \pi+x^{*}} \exp k(x)d y(x),
    \end{split}
\end{equation}
and
\begin{equation}\label{eq181}
\begin{split}
I_2& =\int_{x^{*}}^{2 \pi+x^{*}}(f(x)+ g(x)y(x))\exp\big(2 \ln \frac{y(x)}{y(x^{*})}+k(x)\big)\, dx\\
   & =\frac{1}{(y^*)^2}\int_{x^{*}}^{2 \pi+x^{*}} (f(x)y^2(x)+g(x)y^{3}(x))\exp k(x)\, dx\\
     & =\frac{1}{(y^*)^2}\int_{x^{*}}^{2 \pi+x^{*}} \exp k(x)d y(x).
    \end{split}
\end{equation}
Hence
\begin{equation}\label{eq182}
\begin{split}
G''(y^*) & =-\frac{2}{(y^*)^2}\int_{x^{*}}^{2 \pi+x^{*}} \exp k(x)d y(x)\\
    & =-\frac{2}{(y^*)^2}\big(\int_{x^{*}}^{x_{*}} \exp k(x)d y(x)+\int_{x_{*}}^{2 \pi+x^{*}} \exp k(x)d y(x)\big) \\
   & =-\frac{2}{(y^*)^2}\big(\int_{y^{*}}^{y_{*}} \exp k(x_{1}(y))d y+\int_{y_{*}}^{y^{*}} \exp k(x_{2}(y))d y\big) \\
   & =-\frac{2}{(y^*)^2}\big(\int_{y_{*}}^{y^{*}} (\exp k(x_{2}(y))-\exp k(x_{1}(y)))d y\big),
    \end{split}
\end{equation}
where $x_{1}(y)$ and $x_{2}(y)$ are inverse functions of $y(x)\in[y_{*},y^{*}]$, and $x_{1}(y)\in[x^{*},x_{*}],x_{2}(y)\in[x_{*},x^{*}+2\pi]$. Clearly,
$x_{1}(y)$ and $x_{2}(y)$ are well defined in $[y_{*},y^{*}]$, and $x_{1}'(y)<0$ and $x_{2}'(y)>0$ in $(y_{*},y^{*})$.

Letting $I(y)=k(x_{2}(y))-k(x_{1}(y))$ and $y\in[y_{*},y^{*}]$, it follows from a straightforward computation  that
\begin{equation}\label{eq19}
  I(y)=\int_{x_{1}(y)}^{x_{2}(y)} g(u)y^2(u)\, du,
\end{equation}
and
\begin{equation}\label{eq20}
  I'(y)=g(x_{2}(y))y^{2}(x_{2}(y))x_{2}^{'}(y)-g(x_{1}(y))y^{2}(x_{1}(y))x_{1}^{'}(y).
\end{equation}
It is easy to see that $I(y_{*})=I(y^{*})=0$. Assume that the zero of $g(x_{i}(y))=0$ is $\tilde{y_{i}}$, where $i=1,2$, then we have
$y_{*}<\tilde{y_{2}}, \tilde{y_{1}}<y^{*}$. Without loss of generality, we can assume that $\tilde{y_{2}}\leq \tilde{y_{1}}$. To decide
the sign of $I(y)$, we will divide into three cases to discuss.\\
Case $(i)$: $y_{*}<y\leq \tilde{y_{2}}$. Then we have $g(x_{i}(y))<0, i=1,2$, this leads to $I'(y)<0$, so $I(y)<I(y_{*})=0$. \\
Case $(ii)$: $\tilde{y_{1}}\leq y<y^{*}$. Similar analysis shows that $I'(y)>0$ and $I(y)<I(y^*)=0$.\\
Case $(iii)$: $\tilde{y_{2}}<y<\tilde{y_{1}}$. In this case, we have
\begin{equation}\label{eq21}
\begin{split}
I'(y) =\hspace{1.5mm} & \frac{g(x_{2}(y))}{f(x_{2}(y))+g(x_{2}(y))y(x_{2}(y))}-\frac{g(x_{1}(y))}{f(x_{1}(y))+g(x_{1}(y))y(x_{1}(y))} \\
   =\hspace{1.5mm} & W(y)\bigg(\frac{f(x_{1}(y))}{g(x_{1}(y))}-\frac{f(x_{2}(y))}{g(x_{2}(y))}\bigg),
    \end{split}
\end{equation}
and
\begin{equation}\label{eq22}
\begin{split}
I''(y)=\hspace{1.5mm} & W^{'}(y)\bigg(\frac{f(x_{1}(y))}{g(x_{1}(y))}-\frac{f(x_{2}(y))}{g(x_{2}(y))}\bigg)+W(y)\bigg(\big(\frac{f(x_{1}(y))}{g(x_{1}(y))}\big)^{'}x_{1}^{'}(y)\\
 & -\big(\frac{f(x_{2}(y))}{g(x_{2}(y))}\big)^{'}x_{2}^{'}(y)\bigg),
    \end{split}
\end{equation}
where
$$
W(y)=\frac{g(x_{1}(y))g(x_{2}(y))}{(f(x_{1}(y))+g(x_{1}(y))y(x_{1}(y)))(f(x_{2}(y))+g(x_{2}(y))y(x_{2}(y)))}.
$$

For any $\tilde{y_0}\in  (\tilde{y_{2}}, \tilde{y_{1}})$ satisfying $I'(\tilde{y_0})=0$,  one could find that $W(\tilde{y_0})>0$ and $\frac{f(x_{1}(\tilde{y_0}))}{g(x_{1}(\tilde{y_0}))}-\frac{f(x_{2}(\tilde{y_0}))}{g(x_{2}(\tilde{y_0}))}=0$.
Then, it can be deduced that
$$I''(\tilde{y_0})=W(\tilde{y_0})\bigg(\big(\frac{f(x_{1}(\tilde{y_0}))}{g(x_{1}(\tilde{y_0}))}\big)^{'}x_{1}^{'}(\tilde{y_0})-\big(\frac{f(x_{2}(\tilde{y_0}))}{g(x_{2}(\tilde{y_0}))}\big)^{'}x_{2}^{'}(\tilde{y_0})\bigg)>0.$$
Hence $\tilde{y_0}$ is the minimum point of $I(y)$. In other words, $I(y)$ has only minimum points. On account of $I(\tilde{y_{2}}), I(\tilde{y_{1}})<0$, we have $I(y)\leq \min \{I(\tilde{y_{2}}), I(\tilde{y_{1}})\}<0$.

In sum, $I(y)<0$ for $y\in (y_*, y^*)$, which leads to $G''(y^*)>0$, the limit cycle has the multiplicity of $2$ and is lower stable. The proof of this Lemma is complete. \epf

Equation \eqref{eq4} has the form
$$\frac{dy}{dx}=S(x, y, b),$$
where $S(x, y, b)=(b+c\cos x)y^2+(\sin x-a)y^3$. By means of $\frac{\partial S}{\partial b}=y^2>0$ in the region $y>0$, it is clear that equation \eqref{eq4} is a rotation vector with respect to the parament $b$, so we can study the behaviour of limit cycles when $b$ changes.

\begin{Corollary}\label{cor3.1}
Assume that for some $b_0>0$.  Equation \eqref{eq4} has two hyperbolic limit cycles $y=y_1(x)$ and $y=y_2(x)$, satisfying $y_1(x)<y_2(x)$. If $y_1(x)$ is unstable and $y_2(x)$ is stable, then equation \eqref{eq4} has at least two limit cycles for any $b<b_0$. \end{Corollary}
\noindent\bpf According to Proposition \ref{pro2}, when $b$ decreases, the limit cycle $y_1(x)$ expands and the limit cycle $y_2(x)$ contracts. Thus, if one of these two limit cycles disappears, then it must compose a limit cycle of multiplicity at least two with another limit cycle for some $b=b_*<b_0$. From Lemma \ref{lem3.4}, we know
that this composed limit cycle  is a lower stable limit cycle  of multiplicity two. But when $b$ increases from $b_*$, by Proposition \ref{pro3},
this lower stable limit cycle with multiplicity of two will disappear. This leads to our desires contradiction. Hence  the limit cycles $y_1(x)$ and $y_2(x)$  will exist for any $b<b_0$ and equation \eqref{eq4} has at least two limit cycles.\epf

\begin{Remark}\label{rem3.1} It is easy to check that Lemmas \ref{lem3.2},\ref{lem3.3},\ref{lem3.4} and Corollary \ref{cor3.1}  hold in the region $y<0$.\end{Remark}

Now we consider the case: $0\leq a<1,0\leq b<-c\sqrt{1-a^2},c<0$. By similar discussion, we have the following statements:
\begin{enumerate}
\item[($M_{1}$)] The function $\omega(x)=-\frac{f(x)}{g(x)}$ is monotonically decreasing with respect to $x$;
\item[($M_{2}$)] Lemma \ref{lem3.3} hold;
\item[($M_{3}$)] The multiplicity of the limit cycle of equation \eqref{eq4} is at most two in $y>0$ and $y<0$, respectively, and if it exists, it must be upper stable;
\item[($M_{4}$)] Assume that for some $b_0>0$, equation \eqref{eq4} has two positive (resp. negative) hyperbolic limit cycles $y_1(x)$ and $y_2(x)$, where $y_1(x)<y_2(x)$.
If $y_1(x)$ is stable and $y_2(x)$ is unstable, then equation \eqref{eq4} has at least two positive (resp. negative) limit cycles for any $b>b_0$.
\end{enumerate}

\subsection{ The case: \texorpdfstring{${\xi}g(x)+f(x)$}{}  has no definite sign\label{r}}

 From now on, we will study the case: ${\xi}g(x)+f(x)$  has no definite sign in $[0,2\pi]$, that is, $0<a<1,0<b<|c\sqrt{1-a^2}|$. In fact, we deal with a more general case: $0\leq a<1$ and $0\leq b<|c\sqrt{1-a^2}|$.  Obviously $c\neq 0$, thus we first consider the case $c>0$.

\subsubsection {The case: \texorpdfstring{$c>0$}{} \label{s}}
\noindent \emph{(I) Limit cycles in the region  $y>0$ \label{t}}

For the simplest case $a=b=0$, equation \eqref{eq4} can be written as
 \begin{equation}\label{eq24}
 \frac{dy(x)}{dx}=c\cos x\cdot y^2(x)+\sin x \cdot y^3(x).
\end{equation}
For equation \eqref{eq24}, we have
\begin{Theorem}\label{th5}
Equation \eqref{eq24} has no limit cycle when $c>0$.
\end{Theorem}
\noindent\bpf Suppose that equation \eqref{eq24} has a positive limit cycle $y=y(x)$, then the stability of this limit cycle is determined by
$$II=\int_{0}^{2 \pi} \big(2c\cos x\cdot y(x)+3\sin x\cdot y^{2}(x)\big)\,dx.$$

From \eqref{eq15}, we can get that $\int_{0}^{2 \pi} \big(c\cos x\cdot y(x)+\sin x\cdot y^{2}(x)\big)dx=0$. Hence
\begin{equation}\label{eq25}
\begin{split}
II & =-c \int_{0}^{2 \pi} \cos x \cdot y(x)\,dx\\
  & =-c\big(\int_{\frac{\pi}{2}}^{\frac{3\pi}{2}} \cos x \cdot y(x)\,dx+\int_{\frac{3 \pi}{2}}^{\frac{5\pi}{2}} \cos x \cdot y(x)\,dx\big)\\
  & =-c(\int_{\frac{\pi}{2}}^{\frac{3\pi}{2}} \cos x(y(x)-y(3 \pi-x))\,dx).\\
    \end{split}
\end{equation}
\begin{figure}[!htbp]
\begin{center}
\psfrag{0}{$0$}
\psfrag{1}{$\pi/2$}
 \psfrag{2}{$\pi$}
  \psfrag{3}{$3\pi/2$}
   \psfrag{4}{$2\pi$}
   \psfrag{5}{$5\pi/2$}
   \psfrag{6}{$3\pi$}
 \psfrag{7}{$x$}
 \psfrag{8}{$y$}
  \psfrag{9}{$y(3\pi-x)$}
   \psfrag{a}{$y(x)$}
\psfig{file=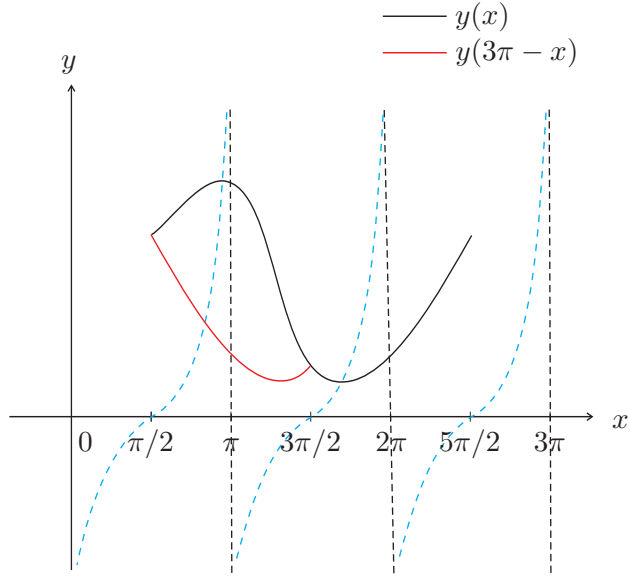,height=3.0in,width=3.2in}
\begin{center}
\caption{Solution curve of equation \eqref{eq24} when $x\in[\frac{\pi}{2},\frac{5\pi}{2}]$.}
\label{Fig. 7}
\end{center}
\end{center}
\end{figure}
Similar to Lemma \ref{lem3.3}, one could obtain the figure of $y(x)$, as shown in Fig.~\ref{Fig. 7}, it is easy to show that $y(x)>y(3 \pi-x)$ for $x\in (\frac{\pi}{2}, \frac{3\pi}{2})$, so $II>0$. That is, all possible limit cycles of equation \eqref{eq24}  must be unstable. However, the zero solution of equation \eqref{eq24} is unstable when $c>0$. Thus, equation \eqref{eq24} has no limit cycle and the theorem is proven. \epf

Next we consider the case: $b=0, 0<a<1$. Now equation \eqref{eq4} becomes
\begin{equation}\label{eq26}
 \frac{dy(x)}{dx}=c\cos x\cdot y^2(x)+(\sin x-a) \cdot y^3(x).
\end{equation}

\begin{Theorem}\label{th3.6}
 Equation \eqref{eq26} has exactly one limit cycle when $0<a<a_{*},c>0$; and no limit cycle  when $a_{*}\leq a<1,c>0$, where $a_{*}$ is the unique zero of $\psi_1(a, c)=0$.
 \end{Theorem}
\noindent\bpf We first analyze the dynamic behavior of system \eqref{eq3} for $0<a<1,b=0,c>0$ to obtain the stability of equation \eqref{eq26} at infinity. We claim that $y=+\infty$ is stable (resp. unstable) for $0<a<a_*$ (resp. $a_*<a<1$), $b=0, c>0$. The complete dynamics of system \eqref{eq3} for $b>0$ is given by Chen and Tang in paper \cite{chen2020global}. In fact, the method in \cite{chen2020global} can be easily extended to the case of $b=0$. By viewing the case $b=0$ as the limit of the case $b>0$, we can obtain the phase portraits for $0<a<a_*$ (resp. $a_*<a<1$), $b=0, c>0$ is $S_{5}$ (resp. $S_{6}$). Then similar to the discussion in Section \ref{infinity}, it is not difficult to conclude that the claim holds.

When $0<a<1,b=0,c>0$, we repeat the proof of Theorem \ref{th5} and obtain that any limit cycle of  equation \eqref{eq26}  must be unstable, as a result, equation \eqref{eq26} has at most one limit cycle. When $0<a<a_{*},c>0$, both $y=0$ and $y=+\infty$ are stable, hence equation \eqref{eq26} has exactly one limit cycle. When $a_{*}<a<1,c>0$, $y=0$ is stable and $y=+\infty$ is unstable, so equation \eqref{eq26} has no limit cycle. In addition, when $a=a_*$, we do not know the stability of $y=+\infty$, by repeating the same idea in the proof of Theorem \ref{th3.1}, it is not difficult to derive that there is also no limit cycle, which finishes the proof.\epf

Essentially, the method of proving Theorem \ref{th5} and  Theorem \ref{th3.6} here is  the same as that in \cite{bravo2009limit}, but for the convenience of readers, we give the detailed proof. At the same time, notice that
$$\frac{\partial \big(c\cos x\cdot y^2(x)+(\sin x-a) \cdot y^3(x)\big)}{\partial a}=-y^3(x)<0,$$
then we have that  equation \eqref{eq26} is a rotation vector with respect to the parameter $a$. Simply discussion shows that when $a$ increases, the limit cycle of equation \eqref{eq26} is born by Hopf bifurcation ($a=0$) and vanishes at $y=+\infty$ ($a=a_*$).

Now we begin to deal with the case $b>0$. The first step is to consider the case: $a<a_{*}$.

\begin{Theorem}\label{th3.7}
Suppose that $0<a<a_{*},0<b<c\sqrt{1-a^{2}},c>0$, there exists a unique $\bar{b}\in(0,c\sqrt{1-a^{2}})$ such that equation \eqref{eq4} has exactly two hyperbolic limit cycles for $0<b<\bar{b}$, one lower stable limit cycle for $b=\bar{b}$, and no limit cycle for $\bar{b}<b<c\sqrt{1-a^{2}}$.
\end{Theorem}
\noindent\bpf When $0<a<a_{*}, b=0, c>0$, the zero solution of equation \eqref{eq4} is stable, and when $b$ increases, the zero solution becomes  unstable, thus
a stable small limit cycle appears; furthermore, when $0<a<a_{*}, b=0, c>0$,  equation \eqref{eq4} has a hyperbolic unstable limit cycle, which cannot disappear as $b$ increases. Consequently, equation \eqref{eq4} has at least two non-zero limit cycles for $0<b\ll1$.

We claim that for $0<a<a_{*}, b>0$ and $c>0$,  equation \eqref{eq4} has at most two limit cycles, counting with multiplicity. Else there exists some $b_0>0$ such that equation \eqref{eq4} has more limit cycles. Without loss of generality, we can assume that these limit cycles are all hyperbolic (else if there exists a non-hyperbolic limit cycle,  then the multiplicity of this non-hyperbolic limit cycle must be two,  we choice $b=b_0-\epsilon$, where $0<\epsilon \ll 1$, this limit cycle will split into two hyperbolic limit cycles).

Suppose that these hyperbolic limit cycles are $y=y_1(x), y_2(x), \ldots, y_m(x)$, $m\geq 3$, where $0<y_1(x)<y_2(x)< \ldots <y_m(x)$. Since the zero solution is unstable,
the limit cycles $y=y_1(x), y_3(x)$ are stable and the limit cycle $y=y_2(x)$ is unstable. Therefore, for $b=b_0$, we find that the unstable  limit cycle $y=y_2(x)$ and the stable limit cycle $y=y_3(x)$ satisfy $y_2(x)<y_3(x)$. It follows from Corollary \ref{cor3.1} that equation \eqref{eq4} has at least two limit cycles for $b<b_0$, but when $b=0$, equation \eqref{eq4} has exactly one limit cycle. This leads to a contradiction.

To sum up, we have shown that if $0<b\ll1$, then equation \eqref{eq4} has exactly two limit cycles $y=y_1(x)$ and $y=y_2(x)$, where $y_1(x)<y_2(x)$, and for generic $b$, then equation \eqref{eq4} has at most two limit cycles.

On the other hand, from Theorem \ref{th3.3}, equation \eqref{eq4} has no limit cycle for $0<a<a_{*},b=c\sqrt{1-a^{2}}$ and $c>0$. Notice that when $b$ increases, the limit cycle $y=y_1(x)$ expands and the limit cycle $y=y_2(x)$ contracts. Hence, there exists a unique $\bar{b}\in(0,c\sqrt{1-a^{2}})$ such that when $b=\bar b$, the limit cycles $y=y_1(x)$ and $y=y_2(x)$ compose a semi-stable limit cycle which is lower stable. Obviously, equation \eqref{eq4} has exactly two hyperbolic limit cycles when $0<b<\bar{b}$; and no limit cycle when $\bar{b}<b<c\sqrt{1-a^{2}}$. The proof of Theorem \ref{th3.7} is finished.\epf

 The last step is to consider the case: $a_{*}\leq a<1$. There is another bifurcation value $\bar a\in (a_{*}, 1)$, which satisfies $\psi_1(\bar a, c)=c\sqrt{1-\bar{a}^2}$.  We have known that $\psi_1(a, c)<c\sqrt{1-a^2}$ for $a_*\leq a<\bar a$ and $c\sqrt{1-a^2}<\psi_1(a, c)$ for $\bar a<a<1$.

 \begin{Theorem}\label{th3.8}
 For equation \eqref{eq4}, we have the following conclusions:
\begin{enumerate}
\item[$(i)$] When $a_{*}\leq a<\bar{a},0<b<\psi_{1}(a,c), c>0$ or $\bar{a}\leq a<1,0<b<c\sqrt{1-a^{2}},c>0$, there is exactly one limit cycle;
\item[$(ii)$] When $a_{*}\leq a<\bar{a}, b=\psi_{1}(a,c), c>0$, there is at most one limit cycle;
\item[$(iii)$] When $a_{*}\leq a< \bar{a}, \psi_{1}(a,c)<b<c\sqrt{1-a^{2}},c>0$, one of the following statements holds:
 \begin{enumerate}
 \item[$(iii.1)$] There is no limit cycle;
 \item[$(iii.2 )$] There exists a unique $\varepsilon \in(0,c\sqrt{1-a^{2}}-\psi_{1}(a,c))$ such that if $\psi_{1}(a,c)<b<\psi_{1}(a,c)+\varepsilon $, then there are exactly two limit cycles; if $b=\psi_{1}(a,c)+\varepsilon $, then there is exactly one lower stable limit cycle; and if $\psi_{1}(a,c)+\varepsilon <b<c\sqrt{1-a^{2}}$, then there is no limit cycle.
 \end{enumerate}
\end{enumerate}
\end{Theorem}
\noindent\bpf (i) When $a_{*}\leq a< \bar{a},0<b<\psi_{1}(a,c),c>0$ or $\bar{a}\leq a<1,0<b<c\sqrt{1-a^{2}},c>0$, both $y=0$  and $y=+\infty$ are unstable, thus equation \eqref{eq4} has at least one limit cycle.

In the following, our goal is to show that equation \eqref{eq4} has at most one limit cycle. Suppose that equation \eqref{eq4} has at least three limit cycles (taking into account multiplicities) for some $b_0>0$. Similar to the proof of Theorem \ref{th3.7}, we will deduce that equation \eqref{eq4} has at least two limit cycles
for $b<b_0$. This is a contradiction with the fact that when $b=0$, so equation \eqref{eq4} has at most one limit cycle.

Specially, for $a_{*}\leq a<\bar{a}$ and $c>0$, the stability of $x=0$ changes when $b$ increases from $b=0$,  a small limit cycle $y=y_b(x)$
appears by Hopf bifurcation, and this limit cycle will expand as $b$ increases.

(ii) Consider the limit function
$$\bar y(x)=\lim_{b\rightarrow \psi_1(a, c)^-}y_b(x).$$

There are two possibilities:

One  possibility is that for some $x_0$, $\bar y(x_0)=+\infty$. That is, the limit cycle  $y=y_b(x)$ is born from $y=0$ and vanishes at $y=+\infty$. Thus, equation \eqref{eq4} has no limit cycle for $b\geq \psi_1(a, c)$ by the no intersection property of rotation vector.

The second  possibility is that $\bar y(x)$ is well defined for each $x$. Then we have obtained a limit cycle $y=\bar y(x)$.

To show the uniqueness of the limit cycle, we first claim that equation \eqref{eq4} has at most two limit cycles (taking into account multiplicities) when $a_{*}\leq a<\bar{a}, b\geq \psi_1(a, c), c>0$.  Else for some $a_0, b_0, c_0$, equation \eqref{eq4} has  at least three limit cycles (we can assume that they are all hyperbolic). From these three limit cycles, we can find one unstable limit cycle $y=y_1(x)$ and one stable limit cycle $y=y_2(x)$, which satisfy $y_1(x)<y_2(x)$. By Corollary \ref{cor3.1}, one could get that equation \eqref{eq4} has at least two limit cycles for $b<b_0$, but when $b=0$, equation \eqref{eq4} has no limit cycle. This leads to a contradiction.

 At last, we claim that equation \eqref{eq4} has  at most one limit cycle. Else we suppose that for some $a_0, b_0=\psi_1(a_0, c_0), c_0$,  equation \eqref{eq4} has  two limit cycles. Let $b$ decrease from $\psi_1(a_0, c_0)$, the two limit cycles will not disappear. This is impossible since equation \eqref{eq4} has one limit cycle when $a=a_0, b<\psi_1(a_0, c_0), c=c_0$.

(iii) From the proof of (ii), we only need to consider the case: equation \eqref{eq4} has exactly one limit cycle when $b=\psi_1(a, c)$. Notice that when $b$ increases from  $\psi_1(a, c)$, the stability of $y=+\infty$ changes, thus one limit cycle (unstable) is bifurcated from $y=+\infty$. That is, when
 $0<b-\psi_{1}(a,c)\ll1$,  equation \eqref{eq4} has at least two limit cycles. We have shown that equation \eqref{eq4} has at most two limit cycles, hence
 equation \eqref{eq4} has exactly two limit cycles $y=y_1(x)$ and $y=y_2(x)$, where $y_1(x)<y_2(x)$. Since $y=0$ is unstable, so $y=y_1(x)$ is stable and $y=y_2(x)$ is unstable. Moreover, as $b$ increase, $y=y_1(x)$ expands and $y=y_2(x)$ contracts.

At last, it follows from Theorem \ref{th3.3} that equation \eqref{eq4} has no limit cycle for $a_{*}\leq a<\bar{a},b=c\sqrt{1-a^{2}}$ and $c>0$. Hence, for equation \eqref{eq4}, there exists a unique $\varepsilon \in(0,c\sqrt{1-a^{2}}-\psi_{1}(a,c))$ such that if $0<b-\psi_{1}(a,c)<\varepsilon $, then there are exactly two limit cycles; if $b-\psi_{1}(a,c)=\varepsilon $, then there is exactly one lower stable limit cycle; and if $\varepsilon <b-\psi_{1}(a,c)<c\sqrt{1-a^{2}}-\psi_{1}(a,c)$, then there is no limit cycle. Thus, we arrive at the conclusion (iii) of the Theorem \ref{th3.8}.

This completes the proof. \epf

\noindent \emph{(II) Limit cycles in the region $y<0$\label{A}}

\begin{Theorem}\label{th3.9}
Equation \eqref{eq26} has no limit cycle in $y<0$ when $0<a<1,c>0$.
\end{Theorem}

The proof is similar to Theorem \ref{th5}, so we omit it here.

\begin{Theorem}\label{th3.10}
Equation \eqref{eq4} has no limit cycle in $y<0$ when $0<a<1,0<b<c\sqrt{1-a^2}$ and $c>0$.
\end{Theorem}
\noindent\bpf Denote the Poincar\'{e} map by $G_b(y_0)=y_b(x_0+2\pi, y_0)-y_0$.  Since when $b=0$, equation \eqref{eq4} has no limit cycle and $y=0$ is stable,
so we have $G_0(y_0)>0$ for all $y_0<0$. Furthermore,  equation \eqref{eq4} is a rotation vector with respect to the parameter $b$, thus one can easily obtain that $\frac{\partial G_b(y_0)}{\partial b}>0$. Immediately, when $b>0$, $G_b(y_0)>G_0(y_0)>0$ for all $y_0<0$, hence equation \eqref{eq4} has no limit cycle in $y<0$, which finishes the proof.\epf

\subsubsection {The case: \texorpdfstring{$c<0$}{} \label{B}}
\noindent \emph{(I) Limit cycles in the region $y>0$ \label{C}}
\begin{Theorem}\label{th3.12}
For equation \eqref{eq4}, the following statements hold:
\begin{enumerate}
\item[$(i)$] If $0<a<a^{*},\psi_{1}(a,c) \leq b<-c\sqrt{1-a^{2}},c<0$, then there is no limit cycle;
\item[$(ii)$] If $0<a<a^{*},0<b<\psi_{1}(a,c),c<0$ or $a^{*}\leq a<1,0<b<-c\sqrt{1-a^{2}},c<0$, then there is exactly one limit cycle.
\end{enumerate}
\end{Theorem}
\noindent\bpf We claim that equation \eqref{eq4} has at most one limit cycle under the above conditions, counting with multiplicity. Else  for some $a_0, b_0, c_0$,  equation \eqref{eq4} has at least two limit cycles.  We can assume that they are all hyperbolic, then we can find two limit cycles $y=y_1(x)$ and $y=y_2(x)$, where $0<y_1(x)<y_2(x)$, such that $y=y_1(x)$ is stable and $y=y_2(x)$ is unstable (note that $y=0$ is unstable). By Statement ($M_{4}$), one could obtain that equation \eqref{eq4} has at least two limit cycles for $a=a_0, b>b_0$ and $c=c_0$, but when $b>-c_0$, equation \eqref{eq4} has no limit cycle. This leads to a contradiction.

$(i)$ When $0<a<a^{*},\psi_{1}(a,c)<b<-c\sqrt{1-a^{2}},c<0$,  $y=0$ is unstable and  $y=+\infty$ is stable, so the number of limit cycles must be even, it is clear that there is no limit cycle. When $b=\psi_{1}(a,c)$, if there is a limit cycle, which must be hyperbolic. Let $b$ increase, this limit cycle still exist, which leads a contradiction.

$(ii)$ When $0<a<a^{*},0<b<\psi_{1}(a,c),c<0$ or $a^{*}\leq a<1,0<b<-c\sqrt{1-a^{2}},c<0$,  both $y=0$ and $y=+\infty$ are unstable,  hence the number of limit cycles must be odd, there is exact one limit cycle.

This completes the proof.\epf

\noindent\emph{(II) Limit cycles in the region $y<0$\label{D}}

\begin{Theorem}\label{th3.15}
For equation \eqref{eq4}, the statements below hold:
\begin{enumerate}
\item[$(i)$] If $0<a<a_{*},\psi_{2}(a,c)\leq b<-c\sqrt{1-a^2},c<0$ or $a_{*}\leq a<1,0<b< -c\sqrt{1-a^2},c<0$, then there is no limit cycle;
\item[$(ii)$] If $0<a<a_{*},0<b<\psi_{2}(a,c),c<0$, then there is exactly one limit cycle.
\end{enumerate}
\end{Theorem}

The proof of Theorem \ref{th3.15} is similar to that of Theorem \ref{th3.12}, thus we omit it.

\section{Proofs of the main results and numerical example\label{Pro}}
In this section, we will give the proof of the main results.

\noindent{\bf Proof of Theorem \ref{mainth2}.}  According to the analysis in Section \ref{existence}, Theorem \ref{mainth2} holds.\epf

\noindent{\bf Proof of Theorem \ref{mainth3}.}  By Theorem \ref{mainth2}, if the Josephson equation has no limit cycle of the first kind, then the sum of the number of limit cycles of the first and second kind is at most two.

If the Josephson equation has one limit cycle of the first kind, then $0\leq a<1, \varphi(a, c)<b<-c\sqrt{1-a^2}$ and $c<0$.
Notice that when $0<a<a_*$ and $c<0$,  $\psi_{2}(a,c)<\varphi(a, c)<\psi_1(a, c)$, we have
$$\{0\leq a<1, \varphi(a, c)<b<-c\sqrt{1-a^2}, c<0\}\subset D_1 \cup D_2,$$
where $D_1=\{0\leq a<a_{*}, \psi_{2}(a,c)<b<-c\sqrt{1-a^{2}}, c<0\}$ and $D_2=\{a_*\leq  a<1, 0<b<-c\sqrt{1-a^{2}}, c<0\}$. It follows from Theorem \ref{th3.12} and \ref{th3.15} that system \eqref{eq4} has at most one limit cycle when $c<0$, thus the sum of the number of limit cycles of the first and second kind is also at most two. Specially, when  $0<a<a^{*},\varphi(a,c)<b<\psi_{1}(a,c),c<0$ or $a^{*}\leq a<1,\varphi(a,c)<b<-c\sqrt{1-a^{2}},c<0$, both kinds of limit cycles exist, thus the upper bound $2$ can be reached.

It is easy to check that if we denote  the configuration of the limit cycle of the Josephson equation by $(i,j)$, where $i$ (resp. $j$) represent the number of limit cycles of the first (resp. second) kind, then the configurations of $(0,2)$ and $(1,1)$ can be realized, which finishes the proof of Theorem \ref{mainth3}.\epf

\begin{figure}[!htbp]
\begin{center}
\psfrag{a}{$A$}
 \psfrag{b}{$B$}
  \psfrag{c}{$C$}
   \psfrag{d}{$L_{1}$}
   \psfrag{e}{$L_{2}$}
\psfig{file=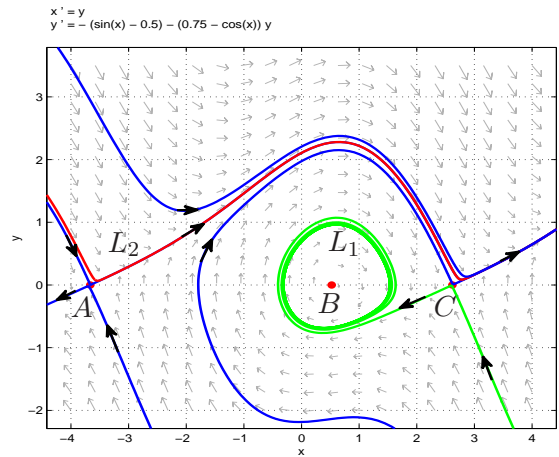,height=3.0in,width=3.2in}
\begin{center}
\caption{Numerical limit cycle of equation \eqref{eq1}.}
\label{Fig. 8}
\end{center}
\end{center}
\end{figure}

In the following, we will show that the Josephson equation has exactly one limit cycle the first kind  and the second kind  by numerical example.
\begin{Example}\label{ex4.1}
Consider $(a,b,c)=(0.5,  0.75,  -1)$. One could obtain that the Josephson equation exactly one first kind limit cycle $L_{1}$ and one second kind limit cycle $L_{2}$, both of which are stable, as shown in Fig.~\ref{Fig. 8}.
\end{Example}


\end{document}